\documentclass[reqno,12pt]{amsart} 

\usepackage[ruled]{algorithm2e}
\renewcommand{\algorithmcfname}{ALGORITHM}
\SetAlFnt{\small}
\SetAlCapFnt{\small}
\SetAlCapNameFnt{\small}
\SetAlCapHSkip{0pt}
\IncMargin{-\parindent}

\usepackage{graphicx}
\usepackage{amssymb,amsmath,amstext}
\usepackage{latexsym}
\usepackage{microtype} 
\usepackage{algorithmic}
\usepackage{xcolor}

\newtheorem{theorem}{Theorem}

\begin{document}

\markboth{J. Brust et al.}{SC-SR1: MATLAB Software for Solving Shape-Changing L-SR1 Trust-Region Subproblems}

\title[Algorithm xxx: SC-SR1: MATLAB Software for Limited-Memory SR1 Trust-Region Methods]
    {Algorithm xxx: SC-SR1: MATLAB Software for Limited-Memory SR1 Trust-Region Methods} 

\author[J. Brust]{Johannes J. Brust}
\email{jbrust@anl.gov}
\address{Argonne National Laboratory}

\author[O. Burdakov]{Oleg Burdakov}
\email{oleg.burdakov@liu.se}
\address{Link\"{o}ping University}

\author[J. Erway]{Jennifer B. Erway}
\email{erwayjb@wfu.edu}
\address{Wake Forest University}

\author[R. Marcia]{Roummel F. Marcia}
\email{rmarcia@ucmerced.edu}
\address{University of California, Merced}

\thanks{
This research is support in part by National Science Foundation grants
CMMI-1333326, CMMI-1334042, IIS-1741264, and IIS-1741490 and
the U.S. Department of Energy, Office of Science, Office of Advanced Scientific Computing Research (ASCR)
under Contract DE-AC02-06CH11347.}


\begin{abstract}
  We present a MATLAB implementation of the 
  symmetric rank-one
  (SC-SR1) method that solves trust-region
  subproblems when a 
limited-memory symmetric rank-one (L-SR1) matrix is used in place of the
true Hessian matrix, which can be used for large-scale optimization. The method takes advantage of two shape-changing
norms~\cite{BYuan02,Burdakov2016} to decompose the trust-region subproblem
into two separate problems. 
  Using one of the
  proposed norms, the resulting subproblems have closed-form
  solutions.   Meanwhile, using the other
  proposed norm, one of the resulting subproblems has a closed-form solution while
  the other is easily solvable
  using techniques that exploit the structure of L-SR1 matrices.
Numerical
  results suggest that the SC-SR1 method is able to solve trust-region
  subproblems to high accuracy even in the so-called ``hard case''.
  {\color{black} When integrated into a trust-region algorithm, extensive numerical experiments suggest that
  the proposed algorithms perform well, when compared with widely used solvers, such as truncated CG.}
\end{abstract}



\maketitle

\makeatletter
\def\BFGS{{\small BFGS}}
\def\LBFGS{{\small L-BFGS}}
\def\LSR{{\small L-SR1}}
\def\SR{{\small SR1}}
\def\CG{{\small CG}}
\def\DFP{{\small DFP}}
\def\LMTR{{\small LMTR}}
\def\OBS{{\small OBS}}
\def\OBSSC{{\small OBS-SC}}
\def\PSB{{\small PSB}}
\def\SCSR{{\small SC-SR1}}
\def\QR{{\small QR}}
\def\MATLAB{{\small MATLAB}}
\newcommand{\minimize}[1]{{\displaystyle\minim_{#1}}}
\newcommand{\minim}{\mathop{\operator@font{minimize}}}
\newcommand{\mgap}{\;\;}
\newcommand{\bgap}{\;\;\;}
\newcommand{\qDef}{{\mathcal Q}}
\newcommand{\defined}{\mathop{\,{\scriptstyle\stackrel{\triangle}{=}}}\,}
\newcommand{\subject}{\mathop{\operator@font{subject\ to}}}  
\newcommand{\words}[1]{\mgap\text{#1}\mgap}
\newcommand\diag{\mathop{\operator@font diag}\nolimits}

\newcommand{\rfm}[1]{\textcolor{black}{#1}}
\newcommand{\obs}[1]{\textcolor{black}{#1}}
\newcommand{\je}[1]{\textcolor{black}{#1}}
\newcommand{\jeo}[1]{\textcolor{black}{#1}}
\newcommand{\jb}[1]{\textcolor{black}{#1}}
\newcommand{\jbo}[1]{\textcolor{black}{#1}}
\renewcommand{\algorithmcfname}{ALGORITHM}

\newcommand{\bp}{\mathbf{p}}
\newcommand{\bv}{\mathbf{v}}
\newcommand{\bfm}[1]{\mathbf{#1}} 
\newcommand{\bk}[1]{\mathbf{#1}_k} 
\newcommand{\bko}[1]{\mathbf{#1}_{k+1}} 
\newcommand{\bsk}[1]{\boldsymbol{#1}_k} 
\newcommand{\bsko}[1]{\boldsymbol{#1}_{k+1}} 
\newcommand{\bs}[1]{\boldsymbol{#1}} 

\newcounter{pseudocode}[section]
\def\thepseudocode{\thesection.\arabic{pseudocode}}
\newenvironment{pseudocode}[2]%
        {%
        \refstepcounter{pseudocode}%
          \AlgBegin %
               {{\bfseries Algorithm \thepseudocode.}\rule[-1.25pt]{0pt}{10pt}#1}%
        #2}%
           {\AlgEnd}

\newcounter{Pseudocode}[section]
\def\thePseudocode{\thesection.\arabic{Pseudocode}}
\newenvironment{Pseudocode}[2]%
        {%
        \refstepcounter{Pseudocode}%
          \AlgBegin %
               {{\bfseries #1.}\rule[-1.25pt]{0pt}{10pt}}%
        #2}%
           {\AlgEnd}

\def\tnu{\tilde{\nu}}

\newcounter{procedureC}
\newcounter{algorithm saved}
\newenvironment{procedureAlg}[1][H]{%
	\setcounter{algorithm saved}{\value{algocf}} 
    	\setcounter{algocf}{\value{procedureC}}
	\renewcommand{\algorithmcfname}{PROCEDURE}
   \begin{algorithm}[#1]%
  }{\end{algorithm}
  \setcounter{procedureC}{\value{algocf}}
  \setcounter{algocf}{\value{algorithm saved}} 
  }

\makeatother

\section{Introduction} \label{intro}

At each iteration of a trust-region
method for minimizing a general nonconvex function $f(\mathbf{x})$, 
the so-called \emph{trust-region subproblem} must
be solved to obtain a step direction:
\begin{equation} \label{eq:trustProblem}
  \minimize{ \mathbf{p}\in\mathbb{R}^n}\mgap\qDef\left(\mathbf{p}\right) \defined \mathbf{g}^T\bp + \frac{1}{2} \bp^T\mathbf{B} \bp
   \bgap\subject \mgap \|\bp\| \le \delta,
\end{equation}
where 
$\mathbf{g}\defined\nabla f\left(\mathbf{x}_k\right)$, $\mathbf{B}$
is an approximation to $\nabla^2 f\left(\mathbf{x}_k\right)$,
$\delta$ is a positive constant, and $\|\cdot\|$ is
a given norm.  
In this article, we describe a \MATLAB{} implementation for solving
the trust-region subproblem (\ref{eq:trustProblem})
when $\bfm{B}$ is a limited-memory symmetric rank-one (\LSR) matrix approximation
of $\nabla^2 f(\bfm{x}_k)$.
In large-scale optimization, solving (\ref{eq:trustProblem}) represents the
bulk of the computational effort in trust-region methods.
The norm used in
(\ref{eq:trustProblem}) not only defines the trust region shape but
also determines the difficulty of solving each subproblem.

The most widely-used norm chosen to define the trust-region subproblem is
the two-norm.  One reason for this choice of norm is that the necessary and
sufficient conditions for a global solution to the subproblem defined by
the two-norm are well-known~\cite{Gay81,MorS83,Sor82}; many methods exploit these
conditions to compute high-accuracy solutions to the trust-region
subproblem (see e.g., \cite{EG10,EGG09,ErwayM14,Gould2010,BruEM15,MorS83}).
The infinity-norm is sometimes 
used to define the subproblem; however, when $\mathbf{B}$ is indefinite, as
can be the case when $\mathbf{B}$ is a \LSR{} matrix,
the subproblem is NP-hard~\cite{MK87,Vav92}.  For more discussion on norms other  than the
infinity-norm we refer the reader to~\cite{ConGT00a}.

\medskip

In this article, we consider the trust-region subproblems defined by 
\emph{shape-changing} norms originally proposed in~\cite{BYuan02}.
Generally speaking, shape-changing norms are norms that depend on
$\mathbf{B}$; thus, in the quasi-Newton setting where the quasi-Newton
matrix $\mathbf{B}$ is updated each iteration, the shape of the trust
region changes each iteration. One of the earliest references to shape-changing norms is 
found in \cite{Gold80} where a norm is implicitly defined by the product of
a permutation matrix and a unit lower triangular matrix that arise from a
symmetric indefinite factorization of $\mathbf{B}$.  Perhaps the most
widely-used shape-changing norm is the so-called ``elliptic norm'' given by
$\|\mathbf{x}\|_A\defined \mathbf{x}^T\mathbf{Ax}$, where $\mathbf{A}$ is a
positive-definite matrix (see, e.g.,~\cite{ConGT00a}).  A well-known use of
this norm is found in the Steihaug method~\cite{Ste83}, 
 and, more
generally, truncated preconditioned conjugate-gradients
(\CG)~\cite{ConGT00a}; these methods
reformulate a two-norm trust-region subproblem
using an elliptic norm to maintain the property that the iterates from
preconditioned \CG{} are increasing in norm.  
Other examples of
shape-changing norms include 
those defined by vectors in the span of $\mathbf{B}$ (see,
e.g.,~\cite{ConGT00a}).

\medskip

The shape-changing norms proposed in~\cite{BYuan02,Burdakov2016} have the
  advantage of breaking the trust-region subproblem into two separate
  subproblems.  Using one of the proposed shape-changing norms, the
  solution of the subproblem then has a closed-form solution.  In the other
  proposed norm, one of the subproblems has a closed-form solution while
  the other is easily solvable.  The publicly-available \LMTR{} codes~\cite{Burdakov2018} solve trust-region subproblems (\ref{eq:trustProblem}) defined by 
these shape-changing norms and the limited-memory Broyden-Fletcher-Goldfarb-Shanno (L-BFGS)
  updates of $\mathbf{B}$.  To our knowledge, there are no other implementations for
  solving trust-region subproblems defined by these shape-changing norms.

\subsection{Overview of the proposed method}
In this paper, we \jb{develop} a \MATLAB{} implementation for solving trust-region \jb{(TR)}
subproblems defined by the two shape-changing norms \jb{described} in~\cite{BYuan02}
when \LSR{} approximations to the Hessian are used 
instead of \LBFGS{} approximations. \jb{For limited-memory algorithms a re-scaling strategy (i.e., effectively re-initializing the Hessian approximation
at each iteration) is often important for the practical performance of the method. Yet, because the structure of \LSR{} matrices
can be exploited to reduce the memory usage even further when a constant initialization is used (i.e., no re-scaling) we provide an option
to chose between such strategies. Moreover, our implementation enables the testing and addition of new solvers by 
swapping out the respective TR subproblem algorithm. In this way, we conduct numerical experiments on large-scale CUTEst 
problems \cite{GouOT03}, comparing the shape-changing methods to truncated CG and an $\ell_2$-norm based algorithm.}
 The proposed method, called the shape-changing SR1 method (\SCSR), \jb{enables high-accuracy subproblem solutions by exploiting the
structure of \LSR{} matrices}.

\medskip

This paper is organized as follows:  In Section 2, we review \LSR{}
matrices, including the compact representation for these matrices and a
method to efficiently compute their eigenvalues and a partial eigenbasis.
In Section 3, we demonstrate how the shape-changing norms decouple the
original trust-region subproblem into two problems and describe the
proposed solver for each subproblem.  Finally, for each shape-changing
norm, we show how to construct a
global solution to (\ref{eq:trustProblem}) from the solutions of
the two decoupled subproblems. Optimality conditions are presented 
for each of these decoupled subproblems in Section 4.   In Section 5,
we demonstrate the accuracy of the proposed solvers, \jb{and compare them on a collection
of large-scale optimization problems.} Concluding
remarks can be found in Section 6.
 
\subsection{Notation}
In this article, the identity matrix of dimension $d$ is denoted by
$\mathbf{I}_d = \left[ \mathbf{e}_1| \cdots | \mathbf{e}_d
\right]$, and depending on the context the subscript $d$ may be
suppressed. Finally, we assume that all \LSR{} updates are computed so that
the \LSR{} matrix is well defined.

\section{L-SR1 matrices}

Suppose $f:\mathbb{R}^n\rightarrow \mathbb{R}$ is a smooth objective function
and $\{ {\bfm{x}}_i\}$, $i=0,\ldots k$, is a sequence of iterates, then the
symmetric rank-one (\SR) matrix is defined using pairs
$(\bfm{s}_i,\bfm{y}_i)$ where
$$
	\mathbf{s}_i\defined \mathbf{x}_{i+1} - \mathbf{x}_i
	\quad 
	\text{and}
	\quad 
	\mathbf{y}_i \defined
	\nabla f(\mathbf{x}_{i+1}) - \nabla f(\mathbf{x}_{i}),
$$ 
and $\nabla f$ denotes the gradient of $f$.
Specifically, given an initial matrix $\mathbf{B}_0$,
$\mathbf{B}_{k+1}$ is defined recursively as
\begin{equation}\label{eq:sr1}
	\mathbf{B}_{k+1} \defined \mathbf{B}_k + \frac{(\mathbf{y}_k- \mathbf{B}_k\mathbf{s}_k)(\mathbf{y}_k - \mathbf{B}_k\mathbf{s}_k)^T}{(\mathbf{y}_k - \mathbf{B}_k\mathbf{s}_k)^T\mathbf{s}_k},
\end{equation}
provided $(\mathbf{y}_k - \mathbf{B}_k\mathbf{s}_k)^T\mathbf{s}_k \ne 0$.
In practice, $\mathbf{B}_0 = \bfm{B}^{(k)}_0$ is often taken to be a scalar multiple of the
identity matrix \jb{that re-scales $\bk{B}$ each iteration}; for the duration of this article we assume that
$\mathbf{B}_0=\gamma_k \mathbf{I}$, $\gamma_k\in\mathbb{R}$.  \emph{Limited-memory} symmetric rank-one
matrices (\LSR) store and make use of only the $m$ most-recently computed
pairs $\{(\mathbf{s}_i,\mathbf{y}_i)\}$, where $m\ll n$ 
(for example, Byrd et al.~\cite{ByrNS94} suggest $m\in
[3,7]$).  For simplicity of notation, we assume that the current
iteration number $k$ is less than the number of allowed stored limited-memory pairs $m$.

The \SR{} update is a member of the Broyden class of updates (see,
e.g.,~\cite{NocW06}).  Unlike widely-used updates such as the
Broyden-Fletcher-Goldfarb-Shanno (\BFGS) and the Davidon-Fletcher-Powell
(\DFP) updates, this update can yield indefinite matrices; that is, \SR{}
matrices can incorporate negative curvature information.  In fact, the \SR{}
update has convergence properties superior to other widely-used
positive-definite quasi-Newton matrices such as \BFGS{};
in particular,~\cite{ConnGT91} give conditions under which
the \SR{} update formula generates a 
sequence of matrices that converge to the true Hessian.
(For more background on the \SR{} update formula, see,
e.g.,~\cite{GrNaS09,KelS98,KhaBS93,NocW06,SunY06,Wol94}.)

\subsection{Compact representation}

The compact representation of \SR{} matrices can be used to
compute the eigenvalues and a partial eigenbasis of these matrices.
In this section, we review the compact formulation of \SR{} matrices.

\medskip

To begin, we define the following matrices:
\begin{eqnarray*}
	\mathbf{S}_k &\defined& [ \ \mathbf{s}_0 \ \ \mathbf{s}_1 \ \ \mathbf{s}_2 \ \ \cdots \ \ \mathbf{s}_{k-1} \ ] \ \in \ \mathbb{R}^{n \times k}, \\
	\mathbf{Y}_k &\defined& [ \ \mathbf{y}_0 \ \ \mathbf{y}_1 \ \ \mathbf{y}_2 \ \ \cdots \ \ \mathbf{y}_{k-1} \ ] \ \in \ \mathbb{R}^{n \times k}.
\end{eqnarray*}
The matrix $\mathbf{S}_k^T\mathbf{Y}_k\in\mathbb{R}^{k \times k}$
can be written as the sum of the following three matrices:
$$
	\mathbf{S}_k^T\mathbf{Y}_k =   \mathbf{L}_k + \mathbf{D}_k + \mathbf{R}_k,
$$
where $\mathbf{L}_k$ is strictly lower triangular, $\mathbf{D}_k$ is diagonal, and $\mathbf{R}_k$ is
strictly upper triangular. 
Then, $\mathbf{B}_{k}$ can be written as 
\begin{equation}\label{eq:form}
	\mathbf{B}_{k} \ = \ \gamma_{k}\mathbf{I} + \mathbf{\Psi}_k \mathbf{M}_k  \mathbf{\Psi}_k^T,
\end{equation}
where $\mathbf{\Psi}_k \in \mathbb{R}^{n \times k}$ and $\mathbf{M}_k \in \mathbb{R}^{k \times
  k}$. 
In particular, $\mathbf{\Psi}_k$ and $\mathbf{M}_k$ are given by
\begin{equation}\label{eqn-PsiM}
	\mathbf{\Psi}_k \ = \ 
       	\mathbf{Y}_k  - \gamma_k\mathbf{S}_k \quad \text{and} \quad
        \mathbf{M}_k \ = \ (\mathbf{D}_k + \mathbf{L}_k + \mathbf{L}_k^T - \gamma_k\mathbf{S}_k^T\mathbf{S}_k)^{-1}.
\end{equation}
The right side of equation (\ref{eq:form}) is the \emph{compact
  representation} of $\mathbf{B}_{k}$; this representation is
due to Byrd et al.~\cite[Theorem 5.1]{ByrNS94}.  For the duration of
this paper, we assume that updates are \jb{made} when both the next
\SR{} matrix $\mathbf{B}_{k}$ is well-defined and $\mathbf{M}_k$
exists \cite[Theorem 5.1]{ByrNS94}.  For notational
simplicity, we assume $\mathbf{\Psi}_k$
has full column rank; when $\mathbf{\Psi}_k$ does not
have full column rank, we refer to~\cite{Burdakov2016} for the 
modifications needed for computing the eigenvalues, which we also review in Section 2.2. 
Notice that the computation of $\mathbf{M}_k$
is computationally admissible since it is a very small symmetric square matrix.

\subsection{\jb{Limited-Memory Updating}}\label{sec-lm}
\jb{For large optimization problems, limited-memory approaches store only a small number of 
vectors to define the \LSR{} representations. Depending on the initialization strategy, specifically
whether $ \gamma_k $ varies between iterations or is constant ($ \gamma_k = \bar{\gamma} $)
the matrices in \eqref{eqn-PsiM} can be effectively stored and updated. We will describe these techniques in subsequent sections.
By setting the parameter $m \ll n$ limited-memory techniques enable inexpensive computations, and replace or insert one column at each
iteration in $ \bfm{Y}_k $ and $ \bfm{S}_k $.
Let an underline below a matrix represent the 
matrix with its first column removed. That is, $ \underline{\bfm{S}}_k $ represents $ \bfm{S}_k $ without its first column.
With this notation, a column update of a matrix, say $ \bk{S} $, by a vector $ \bk{s} $ is defined as follows.
\begin{equation*}
	\text{colUpdate}\left(\bk{S},\bk{s} \right) \defined
	\begin{cases}
		[\: \bk{S} \: \bk{s}\:  ] 						& \text{ if } k < m \\
		[\: \underline{\bfm{S}}_k \: \bk{s}\:  ] 			& \text{ if } k \ge m. \\
	\end{cases}
\end{equation*}
This column update can be implemented efficiently, without copying large amounts of memory,
by appropriately updating a vector that stores index information (``mIdx").  A function to do so is described in
Procedure \ref{alg-colUpdate}:
\begin{procedureAlg}
  \caption{Limited-memory column updating of $ \bk{S} $ by the vector $ \bk{s} $}\label{alg-colUpdate}
\begin{algorithmic}[1]
\ENSURE [$\bk{S}$,\text{mIdx}]=\texttt{colUpdate}($\bk{S}$, $\bk{s}$, mIdx, $m$, $k$);
	\IF {$k = 0$}
	\STATE{$\text{mIdx} \gets \text{zeros}(m,1)$;}	
	\ENDIF
	\IF {$ k < m $}
	\STATE{$ \text{mIdx}(k+1) \gets k+1$;}
	\STATE{$ \bk{S}(:,\text{mIdx}(k+1)) \gets \bk{s} $;}
	\ELSIF {$m \le k$}
	\STATE{$ k_m \gets \text{mIdx}(1) $;}
	\STATE{$ \text{mIdx}(1:(m-1)) \gets \text{mIdx}(2:m) $;}
	\STATE{$ \text{mIdx}(m) \gets k_m $;}
	\STATE{$ \bk{S}(:,\text{mIdx}(m)) \gets \bk{s} $;}
	\ENDIF
\RETURN $\bk{S}$, $\text{mIdx}$;
\end{algorithmic}
\end{procedureAlg} 
Note that this procedure does not copy (or overwrite) large blocks of memory as 
would commands such as $ \{ \bk{S}(:,1:(m-1)) \gets \bk{S}(:,2:m); \bk{S}(:,m) \gets \bk{s} \} $, but instead accesses the relevant 
locations using a stored vector of indices. 
Certain matrix products can also be efficiently updated. As such, the product $ \bk{S}^T \bk{Y} $
does not have to be re-computed from scratch. In order to describe the matrix product updating mechanism,
let an overline above a matrix represent
the matrix with its first row removed. That is, $ \overline{\bfm{S}^{T}_k \bfm{Y}}_k $ represents $ \bfm{S}^{T}_k \bk{Y} $ without its first row.
With this notation, a product update of, say $ \bk{S}^T\bk{Y} $, by matrices $ \bk{S} $, $\bk{Y} $
and vectors $ \bk{s} $, $\bk{y}$ is defined as:
\begin{equation*}
	\text{prodUpdate} \left( \bk{S}^T\bk{Y}, \bk{S}, \bk{Y}, \bk{s}, \bk{y} \right) \defined 
	\begin{cases}
		\left[
			\begin{array}{ c c }
				\bk{S}^T\bk{Y} 		& \bk{S}^T\bk{y} \\
				\bk{s}^T\bk{Y}	& \bk{s}^T \bk{y} 
			\end{array}
		\right] & \text{ if } k < m \vspace{0.1cm} \\		
		\left[
			\begin{array}{ c c }
				\left(\underline{\overline{\bfm{S}^{T}_k \bfm{Y}_k}}\right) 			& 	\underline{\bfm{S}}_k^T\bk{y} \\
				\bk{s}^T \underline{\bfm{Y}}_k					&	 \bk{s}^T \bk{y} 
			\end{array}
		\right] & \text{ if } k \ge m .\\
	\end{cases}
\end{equation*}
This product update can be implemented without recomputing potentially large multiplications, by storing
previous products and information about the column order in $ \bk{S} $ and $ \bk{Y} $. In particular, updating the matrix product is based on storing $ \bk{S}^T \bk{Y} $, $ \bk{S}, \bk{Y} $ and the vector ``mIdx". 
Although a different order is possible, we apply the product update after column updates
of $ \bk{S}, \bk{Y} $ have been done previously. In such a situation the vector, which
stores the appropriate index information (``mIdx") is defined at such a point.}\\

\jb{\begin{procedureAlg}
  \caption{Limited-memory product update $ \bk{S}^T \bk{Y} $ ($ \bk{s} $, $ \bk{y} $ are
  column updates to $ \bk{S} $, $ \bk{Y} $) }\label{alg-prodUpdate}
\begin{algorithmic}[1]
\ENSURE [$\bk{S}^T\bk{Y}$]=\texttt{prodUpdate}($\bk{S}^T\bk{Y}$, $\bk{S}$, $\bk{Y}$, $ \bk{s} $, $ \bk{y} $, mIdx, $m$, $k$);
	\IF {$ k < m $}
	\STATE{$\bk{S}^T\bk{Y}(1:(k+1),k+1) \gets \bk{S}(:,\text{mIdx}(1:(k+1)))^T \bk{y}$;}
	\STATE{$\bk{S}^T\bk{Y}(k+1,1:k) \gets \bk{s}^T\bk{Y}(:,\text{mIdx}(1:k))$;}
	\ELSIF {$m \le k$}
	\STATE{$ \bk{S}^T\bk{Y}(1:(m-1),1:(m-1)) \gets \bk{S}^T\bk{Y}(2:m,2:m) $; \label{alg-prodUpdate:0}}
	\STATE{$ \bk{S}^T\bk{Y}(1:m,m) \gets \bk{S}(:,\text{mIdx}(1:m)^T \bk{y} $; \label{alg-prodUpdate:1}} 
	\STATE{$\bk{S}^T\bk{Y}(m,1:(m-1)) \gets \bk{s}^T\bk{Y}(:,\text{mIdx}(1:(m-1)))$; \label{alg-prodUpdate:2}}
	\ENDIF
\RETURN $\bk{S}^T\bk{Y}$;
\end{algorithmic}
\end{procedureAlg} 
Note that such a product update is computationally much more efficient, than recomputing the
product from scratch. Specifically, when $ m \le k $, the direct product $ \bk{S}^T \bk{Y} $ is
done at $ \mathcal{O}(m^2n) $ multiplications. However, Procedure \ref{alg-prodUpdate} 
does this update with $ \mathcal{O}(2mn) $ multiplications in lines \ref{alg-prodUpdate:1} and \ref{alg-prodUpdate:2},
by reusing previous values from line \ref{alg-prodUpdate:0}. 
Moreover, when the product is symmetric, e.g. 
Procedure \ref{alg-prodUpdate} is invoked by $ \texttt{prodUpdate}( \bk{S}^T\bk{S}, \bk{S}, \bk{S}, \bk{s}, \bk{s}, \text{mIdx}, m, k ) $, 
then $ \bk{S}(:,\text{mIdx}(1:m)^T \bk{s} $ can be stored in line \ref{alg-prodUpdate:1} and reused in line \ref{alg-prodUpdate:2}
(thus only one matrix-vector product is needed, instead of two).
Since limited-memory updating of the \LSR{} matrices varies for the chosen initialization strategy, we describe the 
cases of non-constant initializations $ \gamma_k $ and constant $ \gamma_k = \bar{\gamma} $ next.
\subsubsection{Limited-memory updating of \eqref{eqn-PsiM} using non-constant $ \gamma_k $} \label{sec-up:noncons}
When $ \gamma_k $ varies for every iteration, $ \bsk{\Psi} $ is best implicitly represented by
storing $ \bk{S} $ and $ \bk{Y} $, instead of explicitly forming it (forming $ \bsk{\Psi} $ 
explicitly incurs additional $ \mathcal{O}(mn) $ memory locations 
in $ \bsk{\Psi} = \bk{Y} - \gamma_k \bk{S} $). By
storing the previous $ m $ pairs $ \{ \bfm{s}_i, \bfm{y}_i \}_{i=k-m}^{k-1} $ in the limited-memory 
matrices $ \bk{S} = [\: \bfm{s}_{k-m} \: \cdots \: \bfm{s}_{k-1} \:] \in \mathbb{R}^{n \times m} $ 
and  $ \bk{Y} = [\: \bfm{y}_{k-m} \: \cdots \: \bfm{y}_{k-1} \:] \in \mathbb{R}^{n \times m} $ the 
matrix-vector product  $ \bsk{\Psi}^T \bfm{g} $ (for a vector $\bfm{g}$) is done as
\begin{equation*}
	\bsk{\Psi}^T \bfm{g} = \bk{Y}^T\bfm{g} - \gamma_k ( \bk{S}^T \bfm{g} ).
\end{equation*}
\subsubsection{Limited-memory updating of \eqref{eqn-PsiM} using constant $ \gamma_k = \bar{\gamma} $} \label{sec-up:cons}
When $ \gamma_k = \bar{\gamma} $ is constant, then $ \bk{Y} $ and $ \bk{S} $ do not have to be stored
separately. Instead the limited-memory method stores $m$ previous vectors $ \{ \bs{\psi}_i = \bfm{y}_i - \bar{\gamma} \bfm{s}_i \}_{i=k-m}^{k-1} $,
concatenated in the matrix
\begin{equation*}
	\bsk{\Psi} = \left[ \ \bs{\psi}_{k-m} \ \ \cdots \ \ \bs{\psi}_{k-1} \ \right] \in \mathbb{R}^{n \times m} 
\end{equation*}
Matrix vector products are directly computed as $ \bsk{\Psi}^T \bk{g} $. Subsequently, $ \bk{M} $ from \eqref{eqn-PsiM}
can be updated efficiently by noting that 
\begin{equation*}
	\bk{M}^{-1}\bk{e} = \left(\bk{D} + \bk{L} + \bk{L}^T - \bar{\gamma} \bk{S}^T \bk{S} \right)\bk{e} = \bsk{\Psi}^T\bk{s}.
\end{equation*} 
Because of these simplifications an L-SR1 algorithm with constant initialization strategy can be implemented with about 
half the memory footprint (storing only $ \bsk{\Psi} $ as opposed to $ \bk{Y}, \bk{S} $ (and previous small products)).
However, often the ability to rescale the computations via a non-constant $\gamma_k$ parameter can be advantageous 
in solving large-scale optimization problems. We provide an option to choose between constant or non-constant 
initialization strategies in our implementations.}

\subsection{Eigenvalues}\label{sec-eigs}

In this subsection, we demonstrate how the eigenvalues and a partial
eigenbasis can be computed for \SR{} matrices.  In general, this derivation
can be done for any limited-memory quasi-Newton matrix that admits a
compact representation; in particular, it can be done for any member of the
Broyden convex class~\cite{ErwayM15,B18,BMP18}.  This discussion is based
on~\cite{Burdakov2016}.

Consider the problem of computing the eigenvalues of $\mathbf{B}_{k}$,
which is assumed to be an \LSR{} matrix, obtained from performing $m$
rank-one updates to $\mathbf{B}_0= 
\mathbf{\gamma I}$.  For notational simplicity, we drop subscripts and
consider the compact representation of $\mathbf{B}$:
\begin{equation}\label{eq:Bk} 
	\mathbf{B} = \gamma \mathbf{I} + \mathbf{\Psi} \mathbf{M} \mathbf{\Psi}^T.
\end{equation}
The ``thin'' QR factorization of $\mathbf{\Psi}$ can be written as $\mathbf{\Psi} =
\mathbf{Q}\mathbf{R}$ where $\mathbf{Q} \in \mathbb{R}^{n \times m}$
and $\mathbf{R} \in \mathbb{R}^{m \times m}$ is invertible because,
as it was assumed above, $\mathbf{\Psi}$ has full column rank.  Then,
\begin{equation}\label{eq:eig-1}
	\mathbf{B} = \gamma \mathbf{I} + \mathbf{Q}\mathbf{R}\mathbf{M}\mathbf{R}^T\mathbf{Q}^T.
\end{equation}
The matrix $\mathbf{RMR}^T\in\mathbb{R}^{m\times m}$ is of a
relatively small size, and thus, it is computationally inexpensive to
compute its spectral decomposition.
We define the
spectral decomposition of  $\mathbf{R}\mathbf{M}\mathbf{R}^T$
as $\mathbf{U}\hat{\Lambda}\mathbf{U}^T,$ where $\mathbf{U} \in
\mathbb{R}^{m \times m}$ is an orthogonal
matrix whose columns are made
up of eigenvectors of $\mathbf{R}\mathbf{M}\mathbf{R}^T$ and
$\hat{\Lambda}=\diag(\hat{\lambda}_1,\allowbreak \dots,\allowbreak
\hat{\lambda}_{m})$ is a diagonal matrix
whose entries are the associated eigenvalues.
  
Thus,
\begin{equation}\label{eqn-eig1-nopermute}
	\mathbf{B} = \gamma \mathbf{I} + \mathbf{Q}\mathbf{U}\hat{\Lambda}\mathbf{U}^T\mathbf{Q}^T.
\end{equation}
Since both 
$\mathbf{Q}$ and $\mathbf{U}$ have orthonormal columns, 
$\mathbf{P}_\parallel \defined
\mathbf{Q}\mathbf{U}\in\mathbb{R}^{n\times m}$ also has orthonormal columns.
Let $\bfm{P}_\perp$ denote the matrix whose columns form an orthonormal
basis for $\left({\bfm{P}}_\parallel\right)^\perp$. 
Thus, the spectral decomposition of $\mathbf{B}$ is defined as
	$\mathbf{B} = \mathbf{P}\Lambda_{\gamma} \mathbf{P}^T,$
where
\begin{equation}\label{eqn-PL}
\mathbf{P}\defined \begin{bmatrix}\mathbf{P_\parallel} \,\,\
\mathbf{P_\perp} \end{bmatrix} \quad \text{and}
\quad	\Lambda_{\gamma}  \defined
	\begin{bmatrix}
		\Lambda & 0 \\
		0 & \gamma \mathbf{I}_{n-m}
	\end{bmatrix},
\end{equation}
with $\Lambda_{\gamma}=\diag(\lambda_1,\ldots,\lambda_n)$ and
$\Lambda=\diag(\lambda_1,\ldots,\lambda_{m})=\hat{\Lambda}+\gamma\mathbf{I}\in\mathbb{R}^{m\times m}$. 

\medskip

We emphasize three important properties of the eigendecomposition.
First, all eigenvalues of $\mathbf{B}$ are explicitly obtained and
represented by $ \Lambda_{\gamma} $.  Second, only the first $ m $
eigenvectors of $ \mathbf{B} $ can be explicitly computed, if needed; 
they are represented by $\mathbf{P}_{\parallel} $.  In
particular, since $ \mathbf{\Psi} = \mathbf{Q}\mathbf{R} $, then
	\begin{equation} \label{eqn-matvec-withP}
		\mathbf{P}_{\parallel} = \mathbf{Q}\mathbf{U} = \mathbf{\Psi} \mathbf{R}^{-1}\mathbf{U}.
	\end{equation}
        If $\mathbf{P}_\parallel$ needs to only be available to compute
        matrix-vector products then one can avoid explicitly forming
        $\mathbf{P}_\parallel$ by storing $\mathbf{\Psi}$, $\mathbf{R}$,
        and $\mathbf{U}$. Third, the eigenvalues given by the parameter
          $ \gamma $ can be interpreted as an estimate of the curvature of $ f
          $ in the space spanned by the columns of $ \mathbf{P}_{\perp}$.

          While there is no reason to assume the function $f$ has
          negative curvature throughout the entire subspace
          $\mathbf{P}_\perp$, in this paper, we consider the case
          $\gamma\le 0$ for the sake of completeness.

For the duration of this article, we assume the first $m$ eigenvalues in
$\Lambda_{\gamma}$ are ordered in increasing values, i.e.,
$\Lambda=\diag(\lambda_1,\ldots,\lambda_{m})$ where $\lambda_1\le
\lambda_2 \le \ldots \le \lambda_{m}$ and that $r$ is the multiplicity of $\lambda_1$, i.e.,
$\lambda_1 = \lambda_2 = \cdots = \lambda_r < \lambda_{r+1}$.
For details on updating this partial spectral decomposition
when a new quasi-Newton pair is computed, see~\cite{ErwayM15}.

          \subsection{Implementation}
          
In the above presentation, the {\small QR} factorization was used for \jb{ease of readability} to 
  find a partial spectral decomposition of $\mathbf{B}$. 
  However, there are other approaches that may be better
  suited for different applications.  An alternative approach to computing
  the eigenvalues of $\mathbf{B}$ is presented in~\cite{Lu96} that replaces
  the QR factorization of $\mathbf{\Psi}$ with the SVD and an
  eigendecomposition of a $m\times m$ matrix and $t\times t$
  matrix, respectively, where $t\le m$.  (For more details, 
  see~\cite{Lu96}.)  \jeo{However, experiments in~\cite{BruBEM19} indicate that the QR version of this computation
  outperforms the SVD approach. }
  \jeo{When} $\mathbf{\Psi}^T \bfm{\Psi}$ is
  positive definite (i.e., $\mathbf{\Psi}$ is full rank), the Cholesky
  factorization of $\mathbf{\Psi}^T \bfm{\Psi}=\mathbf{R}^T \bfm{R}$ provides the same $\mathbf{R}$
  needed to form $\mathbf{P}_{\parallel}$ in
  (\ref{eqn-matvec-withP})~\cite{Burdakov2016}.
  \jb{Since $ \bs{\Psi} $ is not explicitly formed when a non-constant initialization $ \gamma = \gamma_k $ is used
  (in this case $ \bs{\Psi} $ is defined by storing $ \bfm{Y}, \bfm{S} $) the product matrix $ \bs{\Psi}^T \bs{\Psi} $ is represented by
  \begin{equation}
  	\label{eq:impsipsi}
  	\bs{\Psi}^T \bs{\Psi} = \bfm{Y}^T \bfm{Y} - 2 \gamma \bfm{Y}^T \bfm{S} + \gamma^2 \bfm{S}^T \bfm{S}
  \end{equation}
  (in \eqref{eq:impsipsi} the matrices $ \bfm{Y}^T\bfm{Y}, \bfm{Y}^T\bfm{S} $ and $ \bfm{S}^T \bfm{S} $
  are stored and updated). In contrast, with a constant initialization $ \gamma_k = \bar{\gamma} $ the product $ \bs{\Psi}^T \bs{\Psi} $
  can be directly updated.}

\medskip

   For the algorithm proposed in this paper,
  it is necessary to be able to compute the eigenvalues of $\bfm{B}$ and to be able to
  compute products with $\bfm{P}_\parallel$.
 However, in our application, it could be the case that  $\mathbf{\Psi}$ is
 not full rank; in
 this case, it is preferable to use
the ${\small \text{LDL}^T}$ decomposition~\cite{gvl96}
  of $\mathbf{\Psi}^T \mathbf{\Psi}$ as 
  proposed in~\cite{Burdakov2016}.  Specifically,
  $$\mathbf{\Pi}^T \mathbf{\Psi}^T\mathbf{\Psi \Pi}=\mathbf{LDL}^T,$$
where $\mathbf{\Pi}$ is a permutation matrix.  If $\mathbf{\Psi}$ is rank-deficient,
    i.e., $\text{rank}(\mathbf{\Psi})=r<m$,
    then
    at least one diagonal entry of $\mathbf{D}$ is zero. \jeo{(In computer arithmetic, it will be relatively small.)}
 In the proposed algorithm, we use the
  following criteria to determine whether entries in $\mathbf{D}$ are sufficiently
  large:
  The $i$th entry of $\mathbf{D}$, i.e.,  $d_i$, is sufficiently large provided
  that
  \begin{equation} \label{eqn-ldlcriteria}
    \jeo{d_i > 10^{-8}\times[\mathbf{\Pi}^T\mathbf{\Psi}^T\mathbf{\Psi\Pi}]_{ii}.}
  \end{equation}
Now, let $J$ to be the set of indices that satisfy (\ref{eqn-ldlcriteria}), \jeo{i.e., $r=|J|$.}
  Furthermore, define $\mathbf{D}_\dagger$ to be the matrix $\mathbf{D}$
  having removed any rows and columns indexed by an element not in $J$
  and 
  $\mathbf{L}_\dagger$ to be the matrix $\mathbf{L}$ having removed columns
  indexed by an element not in $J$.
Then, 
    $$ \mathbf{\Psi}^T \mathbf{\Psi}\approx\mathbf{\Pi L}_\dagger\mathbf{D}_\dagger\mathbf{L}^T_\dagger\mathbf{
      \Pi}^T=\mathbf{\Pi R}_\dagger^T\mathbf{R}_\dagger\mathbf{\Pi}^T ,$$
where $\mathbf{R}_\dagger \defined
  \sqrt{\mathbf{D}_\dagger}\mathbf{L}_\dagger^T\in\mathbb{R}^{r\times m}$.
  Furthermore, 
\begin{equation}\label{eq:eig-1P}
	\mathbf{B} \approx \gamma \mathbf{I} + \mathbf{Q}_\dagger\mathbf{R}_\dagger\mathbf{\Pi}^T\mathbf{M\Pi}\mathbf{R}_\dagger^T\mathbf{Q}_\dagger^T
 \quad \text{with}   \quad  \mathbf{Q}_\dagger\defined
     \left(
       \mathbf{\Psi\Pi}\right)_\dagger\mathbf{R}^{-1}_{\ddagger}\in\mathbb{R}^{n\times
         r},
       \end{equation}\jeo{where}
$\left(\mathbf{\Psi\Pi}\right)_\dagger$ is the matrix  $\mathbf{\Psi\Pi}$ having deleted
any columns indexed by an element not in $J$, and
$\mathbf{R}_\ddagger\in\mathbb{R}^{r\times r}$ is the matrix
$\mathbf{R}_\dagger$ having removed columns indexed by elements not in $J$.
Notice that
the matrix
$\mathbf{R}_\dagger\mathbf{\Pi}^T
\mathbf{M}\mathbf{\Pi}\mathbf{R}_\dagger^T \in\mathbb{R}^{r\times r}$
is full rank.

Thus, the eigenvalue decomposition
$\mathbf{U}\hat{\Lambda}\mathbf{U}^T$ is now computed not for $\mathbf{RMR}^T$
as in Section \ref{sec-eigs}, but for
$\mathbf{R_{\dagger}\Pi}^T \mathbf{M\Pi R_{\dagger}}^T$.
Furthermore,  $\bfm{P}_\parallel$ in (\ref{eqn-matvec-withP})
is computed as
\begin{equation}\label{eqn-Pparallel}
  \mathbf{P}_\parallel = \mathbf{Q_\dagger U} = \left(\mathbf{\Psi\Pi}\right)_\dagger \mathbf{R}_{\ddagger}^{-1} \mathbf{U}
\end{equation}
\jb{when a constant initialization is used (since $ \bs{\Psi} $ is explicitly formed), and as
\begin{equation}\label{eqn-Pparallel_nc}
	\mathbf{P}_\parallel = \mathbf{Q_\dagger U} = \left(\bfm{Y\Pi}\right)_\dagger \mathbf{R}_{\ddagger}^{-1} \mathbf{U} 
				- \gamma\left(\bfm{S\Pi}\right)_\dagger \mathbf{R}_{\ddagger}^{-1} \mathbf{U}
\end{equation}
when a non-constant initialization is used.}

\medskip
 Algorithm 1
 details the computation of the elements needed to form $\mathbf{P}_\parallel$,
 using the ${\small \text{LDL}^T}$ decomposition.  It produces
 $\mathbf{\Lambda}$, $\mathbf{R}_{\ddagger}$, $\mathbf{U}$, and $\mathbf{\Pi}$.
There are several pre-processing and post-processing steps in this algorithm.
Namely, lines 7 and 9 are used to remove any spurious complex round-off
error,
line 10 is to order the eigenvalues and associated eigenvectors, and line
12
sets any small eigenvalue (in absolute value) to zero.  An alternative 
to forming and storing $\bfm{R}_{\ddagger}$ is to maintain
 $\bfm{R}_{\dagger}$ and the index set $J$. \jb{Moreover, since it is typically more efficient
 to update the product $ \bs{\Psi}^T \bs{\Psi} $ instead of forming it from scratch, the argument ``$ \Psi \lor \Psi^T \Psi $"
 is used to enable passing either of the two inputs, depending on the context.}
\begin{algorithm}[H]
  \caption{Computing $\mathbf{R}_\ddagger$, $\mathbf{\Lambda}$,
      $\mathbf{U}$, and  $\mathbf{\Pi}$ using the ${\small \text{LDL}^T}$ decomposition}\label{alg-ldl}
\begin{algorithmic}[1]
\ENSURE [$\mathbf{R}_\ddagger$,
$\mathbf{\Lambda}$, $\mathbf{U}$, $\mathbf{\Pi}$, \jeo{$J$}]=\texttt{ComputeSpectral}(\jb{$\Psi \lor \Psi^T \Psi $}, $\mathbf{M}^{-1}$, $\gamma$, $\tau$);
  \STATE Compute the ${\small \text{LDL}^T}$ decomposition of $\bs{\Psi}^T
  \bs{\Psi}$ and store the factors $\mathbf{L}$ and $\mathbf{D}$ matrices, and store
$\mathbf{\Pi}$ (as a vector with the permutation information);
\STATE
   Find the indices of elements of $\mathbf{D}$ that are sufficiently large  using
   (\ref{eqn-ldlcriteria})  and store as
   $J$;
\STATE Form $\mathbf{D}_\dagger$ by storing the rows and columns of $\mathbf{D}$ corresponding to
indices of $J$;
\STATE   Form $\mathbf{L}_\dagger$ by storing the columns of $\mathbf{L}$
corresponding the indices of $J$;
\STATE
   $\mathbf{R}_\dagger \gets \sqrt{\mathbf{D}_\dagger}\mathbf{L}_\dagger^T$; 
\STATE $\mathbf{T}\gets  \mathbf{R}_\dagger\mathbf{\Pi}^T\mathbf{M\Pi R}^T_\dagger$;
   \STATE   Compute the spectral decomposition $\mathbf{U}\mathbf{\hat{\Lambda}}\mathbf{U}^T$of $(\mathbf{T}+\mathbf{T}^T)/2$;
\STATE Form $\mathbf{R}_\ddagger  $ by storing the columns of $\mathbf{R}_\dagger$
corresponding
to columns of $J$;
\STATE $\mathbf{\hat{\Lambda}}\gets \text{real}(\mathbf{\hat{\Lambda}})$
\STATE Order the entries in $\mathbf{\hat{\Lambda}}$ from low to high and
rearrange  
the columns of $\mathbf{U}$ accordingly to maintain the spectral
decomposition
of $(\mathbf{T}+\mathbf{T}^T)/2$;
\STATE $\mathbf{\Lambda}\gets \mathbf{\hat{\Lambda}}+\gamma \bfm{I}$;
\IF {$|\mathbf{\Lambda}_{ii}|<\tau$ for any $i$ } \STATE 
{$\mathbf{\Lambda}_{ii}\gets 0$};
\ENDIF
\RETURN $\mathbf{R}_\ddagger$, $\mathbf{\Lambda}$, $\mathbf{U}$,  $\mathbf{\Pi}$;
\end{algorithmic}
\end{algorithm}

The output of Algorithm 1 includes the factors
 of $\mathbf{P}_\parallel$ (see
(\ref{eqn-Pparallel})), i.e.,
  $\mathbf{R}_\ddagger$, $\mathbf{U}$, and
  $\mathbf{\Pi}$, \jeo{as well as $J$.}
  For the method proposed in Section 3, products 
with $\mathbf{P}_\parallel$ are computed as a sequence of explicit
matrix-vector products with the factors of $\mathbf{P}_\parallel$.
In practice, the permutation matrix $\mathbf{\Pi}$ is not stored explicitly;
instead, the permutation is applied implicitly using a vector that
maintains
the order of the columns after the permutation matrix is applied.  Thus,
products with $\mathbf{P}_\parallel$ are computed using only matrix-vector
products together with a rearranging of columns.



\section{Proposed method} \label{sec-proposed}
The proposed method is able to solve the \LSR\ trust-region subproblem to
high accuracy, even when $\mathbf{B}$ is indefinite.  The method makes use
of the eigenvalues of $\mathbf{B}$ and the factors of
$\mathbf{P}_\parallel$.  To describe the method, we first transform the
trust-region subproblem \eqref{eq:trustProblem} so that the quadratic
objective function becomes separable.  Then, we describe the shape-changing
norms proposed in \cite{BYuan02,Burdakov2016} that decouples the
separable problem into two minimization problems, one of which has a
closed-form solution while the other can be solved very efficiently.
Finally, we show how these solutions can be used to construct a solution to
the original trust-region subproblem.

\subsection{Transforming the Trust-Region Subproblem}

Let $\mathbf{B} = \mathbf{P}\Lambda_{\gamma}\mathbf{P}^T$ be the
eigendecomposition of $\mathbf{B}$ described in Section 2.2.  Letting
$\mathbf{v} = \mathbf{P}^T\mathbf{p}$ and $\mathbf{g}_{\mathbf{P}} =
\mathbf{P}^T\mathbf{g}$, the objective function
$\mathcal{Q}(\mathbf{p})$ in \eqref{eq:trustProblem} can be written as
a function of $\bfm{v}$:
\begin{equation*}
		\qDef\left( \mathbf{p} \right) 
		= \mathbf{g}^T\mathbf{p} + \frac{1}{2} \mathbf{p}^T \mathbf{B} \mathbf{p} 
		= \mathbf{g}^T_{\mathbf{P}}\mathbf{v} + \frac{1}{2}\mathbf{v}^T \Lambda_{\gamma} \mathbf{v} \defined q\left( \mathbf{v} \right).
\end{equation*}
With $ \mathbf{P} = \left[ \mathbf{P}_{\parallel} \quad \mathbf{P}_{\perp} \right] $, we partition $\mathbf{v}$ and $\mathbf{g}_{\mathbf{P}}$ as follows:
\begin{equation*}
	\mathbf{v} = \mathbf{P}^T\mathbf{p} = 
	\left[
		\begin{array}{c}
			\mathbf{P}^T_{\parallel} \mathbf{p} \\
			\mathbf{P}^T_{\perp} \mathbf{p} \\
		\end{array}
	\right]
	=
	\left[
		\begin{array}{c}
			\mathbf{v}_{\parallel}\\
			\mathbf{v}_{\perp} \\
		\end{array}
	\right]
	\quad \text{and} \quad
	\mathbf{g}_{\mathbf{P}} = 
		\left[ \begin{array}{c} 
				\mathbf{P}^T_{\parallel} \mathbf{g} \\
				\mathbf{P}^T_{\perp}\mathbf{g}
			\end{array} 
		\right] =
		\left[ \begin{array}{c} 
			\mathbf{g}_{\parallel} \\
			\mathbf{g}_{\perp}
		\end{array} \right],
\end{equation*}
where $\mathbf{v}_{\parallel}, \mathbf{g}_{\parallel} \in \mathbb{R}^{m}$ and $\mathbf{v}_{\perp}, \mathbf{g}_{\perp} \in \mathbb{R}^{n-m}$.
Then,
\begin{eqnarray}		\label{eq:q_decomp}
	q\left( \mathbf{v} \right)  
		&=&
		\begin{bmatrix}
			\mathbf{g}_{\parallel}^T \ & \mathbf{g}_{\perp}^T
		\end{bmatrix}
		\begin{bmatrix}
			\mathbf{v}_{\parallel} \\ \mathbf{v}_{\perp}
		\end{bmatrix}
		+ 
		\frac{1}{2}
		\begin{bmatrix}
			\mathbf{v}_{\parallel}^T \  & \mathbf{v}_{\perp}^T
		\end{bmatrix}
		\begin{bmatrix}
			\Lambda \\ 
			& \gamma\mathbf{I}_{n-m} 
		\end{bmatrix}
		\begin{bmatrix}
			\mathbf{v}_{\parallel} \\ \mathbf{v}_{\perp}
		\end{bmatrix}			
		\nonumber \\
		&=& \mathbf{g}^T_{\parallel}\mathbf{v}_{\parallel} + \mathbf{g}^T_{\perp}\mathbf{v}_{\perp} 
				+ \frac{1}{2} \left( \mathbf{v}^T_{\parallel}\Lambda\mathbf{v}_{\parallel} 
				+ \gamma \left\| \mathbf{v}_{\perp} \right\|^2 \right) \nonumber \\
		&=& q_{\parallel} \left( \mathbf{v}_{\parallel} \right) + q_{\perp} \left( \mathbf{v}_{\perp} \right),
\end{eqnarray}
where
$$
		q_{\parallel} \left( \mathbf{v}_{\parallel} \right) 
		\defined 
		\mathbf{g}^T_{\parallel}\mathbf{v}_{\parallel} + \frac{1}{2}\mathbf{v}^T_{\parallel}\Lambda\mathbf{v}_{\parallel}
		\quad \text{and} \quad
		q_{\perp} \left( \mathbf{v}_{\perp} \right) 
		\defined 
		\mathbf{g}^T_{\perp}\mathbf{v}_{\perp} + \frac{\gamma}{2} \left\| \mathbf{v}_{\perp} \right\|^2.
$$
Thus, the trust-region subproblem \eqref{eq:trustProblem} can be expressed as
\begin{equation}	\label{eq:qv}
		\minimize{\| \mathbf{ Pv } \| \le \delta }  \: \: 
		q\left( \mathbf{v}\right) 
		=  
		\left\{ q_{\parallel}\left( \mathbf{v}_{\parallel}\right) + q_{\perp}\left( \mathbf{v}_{\perp}\right) \right\}.
\end{equation}
Note that the function $q(\mathbf{v})$ is now separable in
$\mathbf{v}_{\parallel}$ and $\mathbf{v}_{\perp}$.  To completely decouple
\eqref{eq:qv} into two minimization problems, we use a shape-changing norm
so that the norm constraint $\| \mathbf{ Pv } \| \le \delta$ decouples into
separate constraints, one involving $\mathbf{v}_{\parallel}$ and the other
involving $\mathbf{v}_{\perp}$.

\subsection{Shape-Changing Norms}
Consider the following shape-changing norms proposed in~\cite{BYuan02,Burdakov2016}:
 \begin{align}
	\label{eq:sc_2}
	\| \mathbf{p} \|_{\mathbf{P},2} 
	&\defined
	\max\left( \| \mathbf{P}_{\parallel}^T\mathbf{p} \|_2, \| \mathbf{P}_{\perp}^T\mathbf{p} \|_2 \right )
	\hspace{.125cm} = \max\left(  \| \mathbf{v}_{\parallel} \|_2, \| \mathbf{v}_{\perp} \|_2\right), \\
	\label{eq:sc_inf}
	\| \mathbf{p} \|_{\mathbf{P},\infty} 
	&\defined
	\max\left( \| \mathbf{P}_{\parallel}^T\mathbf{p} \|_{\infty}, \| \mathbf{P}_{\perp}^T\mathbf{p} \|_2 \right )
	= \max\left(  \| \mathbf{v}_{\parallel} \|_{\infty} , \| \mathbf{v}_{\perp} \|_2 \right).
 \end{align}
We refer to them as the $(\mathbf{P},2)$ and the $(\mathbf{P},\infty)$ norms, respectively.
Since $\mathbf{p} = \mathbf{Pv}$, the trust-region constraint in \eqref{eq:qv} can be expressed in these norms as
\begin{eqnarray*}
	\| \mathbf{Pv} \|_{\mathbf{P},2} \ \le \delta \quad &\text{if and only if}& \quad
	\| \mathbf{v}_{\parallel} \|_2 \ \le \delta \ \text{and} \
	\| \mathbf{v}_{\perp} \|_2 \le \delta, \\
	\| \mathbf{Pv} \|_{\mathbf{P},\infty}  \le \delta \quad &\text{if and only if}& \quad
	\| \mathbf{v}_{\parallel} \|_{\infty} \le \delta \ \text{and} \
	\| \mathbf{v}_{\perp} \|_2 \le \delta.
\end{eqnarray*}
Thus, from \eqref{eq:qv}, the trust-region subproblem is given
for the $(\mathbf{P},2)$ norm by
\begin{equation} \label{eq:qv2}
	\minimize{\| \mathbf{ Pv } \|_{\mathbf{P},2} \le \delta }  \: \: 
	q\left( \mathbf{v}\right) 
	=  
	\minimize{\| \mathbf{ v }_{\parallel} \|_{2} \le \delta }  \: \: 
	q_{\parallel}\left( \mathbf{v}_{\parallel}\right) 
	+ 
	\minimize{\| \mathbf{ v }_{\perp} \|_{2} \le \delta }  \: \: 
	q_{\perp}\left( \mathbf{v}_{\perp}\right),
\end{equation}
and using the $(\mathbf{P},\infty)$ norm it is given by
\begin{equation} \label{eq:qvinfty}
	\minimize{\| \mathbf{ Pv } \|_{\mathbf{P},\infty} \le \delta }  \: \: 
	q\left( \mathbf{v}\right) 
	=  
	\minimize{\| \mathbf{ v }_{\parallel} \|_{\infty} \le \delta }  \: \: 
	q_{\parallel}\left( \mathbf{v}_{\parallel}\right) 
	+ 
	\minimize{\| \mathbf{ v }_{\perp} \|_{2} \le \delta }  \: \: 
	q_{\perp}\left( \mathbf{v}_{\perp}\right).
\end{equation} 
As shown in~\cite{Burdakov2016}, these norms are equivalent to the 
two-norm, i.e.,
\begin{eqnarray*}
	\frac{1}{\sqrt{2}} \| \mathbf{p} \|_2 \  \le  & \| \mathbf{p} \|_{\mathbf{P},2}  & \le \  \| \mathbf{p} \|_2 \\
	\frac{1}{\sqrt{m}} \| \mathbf{p} \|_2 \ \le &  \| \mathbf{p} \|_{\mathbf{P},\infty} & \le \ \| \mathbf{p} \|_2.
\end{eqnarray*}
Note that the latter equivalence factor depends on the number of
stored quasi-Newton pairs $m$ and not on the number of variables ($n$).

Notice that the shape-changing norms do not place equal
  value on the two subspaces since the region defined by the subspaces is of
  different size and shape in each of them. However, because of norm
  equivalence, the shape-changing region insignificantly differs from the
  region defined by the two-norm, the most commonly-used choice of norm.  \medskip

We now show how to solve the decoupled subproblems.

\subsection{Solving for the optimal $\mathbf{v}_{\perp}^*$}
\label{sec:v_perp}

The subproblem
\begin{equation}
	\label{eq:sp_vperp_2}
	\minimize{\| \mathbf{v}_{\perp} \|_2 \le \delta} \quad q_{\perp}\left ( \mathbf{v}_{\perp}\right) \equiv 
	\mathbf{g}^T_{\perp}\mathbf {v}_{\perp} + \frac{\gamma}{2} \| \mathbf{v}_{\perp} \|_2^2 
\end{equation}
appears in both \eqref{eq:qv2} and \eqref{eq:qvinfty}; its optimal solution
can be computed by formula.  For the quadratic subproblem
(\ref{eq:sp_vperp_2}) the solution $\mathbf{v}_{\perp}^*$ must satisfy the
following optimality conditions found in~\cite{Gay81,MorS83,Sor82} associated
with \eqref{eq:sp_vperp_2}:
For some $\sigma_{\perp}^* \in
\mathbb{R}^+$,
\begin{subequations}
	\label{eq:opt_vperp}
	\begin{align}
		\label{eq:vperp_c1}
	  		\left( \gamma+\sigma^{*}_{\perp} \right)\mathbf{v}^*_{\perp} 							
	  		&= -\mathbf{g}_{\perp}, \\
	  	\label{eq:vperp_c2}
	  		\sigma^{*}_{\perp} \left( \|\mathbf{v}^*_{\perp}\|_2 - \delta \right) 			
			&= 0, \\
	  	\label{eq:vperp_c4}
	  		\|\mathbf{v}^*_{\perp }\|_2
			& \le \delta, \\
	  	\label{eq:vperp_c5}
	  		\gamma + \sigma^{*}_{\perp}  
			& \ge 0.
	\end{align}
\end{subequations}	

Note that the optimality conditions are satisfied by
$(\mathbf{v}_{\perp}^*, \sigma_{\perp}^*)$ given by
\begin{equation}
	\label{eq:soln_vperp}
	\mathbf{v}^*_{\perp} = 
	\begin{cases}			
		-\frac{1}{\gamma} \mathbf{g}_{\perp}
		& \text{ if } \gamma > 0 \text{ and } \left \| \mathbf{g}_{\perp} \right \|_2 \le \delta | \gamma|, \\
		\delta \mathbf{u} 											
		& \text{ if } \gamma \le 0 \text{ and } \| \mathbf{g}_{\perp} \|_2 = 0, \\
		-\frac{ \delta}{\| \mathbf{g}_{\perp} \|_2}    \mathbf{g}_{\perp} 
		& \text{ otherwise, }			
	\end{cases}
\end{equation}
and
\begin{equation} 
	\label{eq:soln_sigmaperp}
	\sigma^*_{\perp} = 
	\begin{cases}			
		0 									
		& \text{ if } \gamma > 0 \text{ and }  \left \| \mathbf{g}_{\perp} \right \|_2 \le \delta |\gamma|, \\
		\frac{\left\| \mathbf{g}_{\perp} \right\|_2}{\delta} - \gamma   
		& \text{ otherwise,}			
	\end{cases}
\end{equation}
where $\mathbf{u} \in \mathbb{R}^{n-m}$ is any unit vector
with respect to the two-norm.

\subsection{Solving for the optimal $\mathbf{v}_{\parallel}^*$}	
\label{sec:v_parallel}

\noindent In this section, we detail how to solve for the
optimal $\mathbf{v}_\parallel^*$ when either the $(\mathbf{P}, \infty)$-norm
or the $(\mathbf{P},2)$-norm is used to define the trust-region subproblem.

\bigskip

\noindent 
\textbf{$(\mathbf{P},\infty)$-norm solution}. 
If the shape-changing $(\mathbf{P}, \infty)$-norm is used in \eqref{eq:qv}, then the subproblem in $ \mathbf{v}_{\parallel} $ is
\begin{equation}
	\label{eq:sp_vpar_inf}
	\minimize{\left\| \mathbf{v}_{\parallel} \right\|_{\infty} \le \delta} \quad q_{\parallel}\left( \mathbf{v}_{\parallel} \right) 
	= 
	\mathbf{g}^T_{\parallel} \mathbf{v}_{\parallel} 
	+ 
	\frac{1}{2} \mathbf{v}^T_{\parallel} \Lambda \mathbf{v}_{\parallel}.
\end{equation}
The solution to this problem is computed by separately minimizing $ m $ scalar quadratic problems of the form
\begin{equation}
	\label{eq:sp_vpar_infi}
	\minimize{ |[ \mathbf{v}_{\parallel}]_i|\le \delta } 
	\quad  
	q_{\parallel,i}( [\mathbf{v}_{\parallel}]_i) = 
	\left[ \mathbf{g}_{\parallel} \right]_i\left[ \mathbf{v}_{\parallel} \right]_i 
	+ 
	\frac{\lambda_i}{2} \left( \left[ \mathbf{v}_{\parallel} \right]_i \right)^2, 
	\quad \quad 1 \le i \le m.
\end{equation}
The minimizer depends on the convexity of $q_{\parallel,i}$, i.e., the sign of $\lambda_i$.  The solution to \eqref{eq:sp_vpar_infi} is given as follows: 
\begin{equation}
	\label{eq:vpar_inf}
	[\mathbf{v}^*_{||} ]_i =
	\begin{cases}
		-\frac{\left[ \mathbf{g}_{||}\right]_i}{ \lambda_i} 
		& \text{ if } \left| \frac{ \left[ \mathbf{g}_{||}\right]_i  }{\lambda_i} \right| \le \delta \text{ and } \lambda_i > 0, \\
		c
		& \text{ if } \left[ \mathbf{g}_{\parallel}\right]_i = 0, \: \lambda_i = 0,\\
		- \text{sgn}(\left [ \mathbf{g}_{\parallel} \right ]_i) \delta											
		& \text{ if } \left[ \mathbf{g}_{\parallel}\right]_i \ne 0, \: \lambda_i = 0,\\				
		\pm\delta 											
		& \text{ if } \left[ \mathbf{g}_{\parallel}\right]_i = 0, \: \lambda_i < 0,\\
		-\frac{\delta}{ \left| \left[ \mathbf{g}_{||}\right]_i \right|} \left[ \mathbf{g}_{||}\right]_i 
		& \text{ otherwise},
	\end{cases}
\end{equation}
where $c$ is any real number in $[-\delta, \delta]$ and ``sgn'' denotes the signum
function (see~\cite{Burdakov2016} for
details).
	
\bigskip	

\noindent \textbf{$(\mathbf{P},2)$-norm solution}: 
If the shape-changing $(\mathbf{P}, 2)$-norm is used in \eqref{eq:qv}, then the subproblem in $ \mathbf{v}_{\parallel} $ is
\begin{equation}
		\label{eq:sp_vpar_2}
		\minimize{\| \mathbf{v}_{\parallel} \|_2 \le \delta} \quad 
		q_{\parallel}\left( \mathbf{v}_{\parallel} \right) = 
		\mathbf{g}^T_{\parallel} \mathbf{v}_{\parallel} 
		+ 
		\frac{1}{2} \mathbf{v}^T_{\parallel} \Lambda \mathbf{v}_{\parallel}.		
\end{equation}
The solution $\mathbf{v}^*_{\parallel}$ must satisfy the following optimality conditions~\cite{Gay81,MorS83,Sor82} associated with (\ref{eq:sp_vpar_2}):
For some $\sigma_\parallel^*\in\mathbb{R}^+$,

\begin{subequations}
	\label{eq:opt_vpar}
	\begin{align}
		\label{eq:vpar_c1}
	  	(\Lambda+\sigma^{*}_{\parallel}\mathbf{I})\mathbf{v}^*_{\parallel} 						
		&= -\mathbf{g}_{\parallel}, \\
	  	\label{eq:vpar_c2}
	  	\sigma^{*}_{\parallel} \left( \|\mathbf{v}^*_{\parallel}\|_2 - \delta \right) 			
		&= 0, \\
	  	\label{eq:vpar_c4}
	  	\|\mathbf{v}^*_{\parallel}\|_2														
		& \le \delta, \\
	  	\label{eq:vpar_c5}
	  	\lambda_i + \sigma^{*}_{\parallel}  											
		& \ge 0 \quad \text{ for } 1 \le i \le m.
	\end{align}
\end{subequations}
A solution to the optimality conditions
(\ref{eq:vpar_c1})-(\ref{eq:vpar_c5}) can be computed using the method
found in~\cite{BruEM15}.  For completeness, we outline the method here;
this method depends on the sign of $\lambda_1$. 
Throughout these cases, we make use of the 
expression of $\mathbf{v}_{\parallel}$ as a function of $\sigma_\parallel$.
That is, from
the first optimality condition \eqref{eq:vpar_c1}, we write
\begin{equation}
	\label{eq:v_sigma}
	\mathbf{v}_{\parallel}\left( \sigma_{\parallel} \right) 
	= 
	-\left( \Lambda + \sigma_{\parallel} \mathbf{I} \right)^{-1} \mathbf{g}_{\parallel},
\end{equation}
with $\sigma_{\parallel} \ne -\lambda_i$ for $1 \le i \le m$.

\bigskip
	
\noindent 
\textbf{Case 1 ($\lambda_1 > 0 $).}  When $\lambda_1>0$,
the unconstrained minimizer is computed (setting $\sigma_\parallel^*=0$):
\begin{equation}
	\label{eq:v_par_solpos}
	\mathbf{v}_{\parallel} \left( 0 \right) 
	= 
	- \Lambda^{-1}\mathbf{g}_{\parallel}.
\end{equation}
 If $\mathbf{v}_{\parallel}(0)$ is feasible, 
i.e., $\| \mathbf{v}_{\parallel} \left( 0 \right)  \|_2 \le \delta $ 
then $\mathbf{v}_{\parallel}^* = \mathbf{v}_{\parallel} (0)$ is the global
minimizer; otherwise, $ \sigma^*_{\parallel} $ is the solution to the secular
equation \eqref{eq:secular} (discussed below).  The minimizer to the
problem \eqref{eq:sp_vpar_2} is then given by
\begin{equation}
		\label{eq:v_secularsoln}
		\mathbf{v}^*_{\parallel} = -\left( \Lambda + \sigma^*_{\parallel}\mathbf{I}\right)^{-1}\mathbf{g}_{\parallel}.
\end{equation}

\bigskip

\noindent \textbf{Case 2 ($ \lambda_1 = 0 $).} 
If $ \mathbf{g}_{\parallel} $ is in the range of $\Lambda$, i.e., $[\mathbf{g}_{\parallel}]_i = 0$ for $ 1 \le i \le r$, then set $\sigma_{\parallel} = 0$ and let 
\begin{equation*}
			\mathbf{v}_{\parallel}\left( 0 \right) = - \Lambda^{\dagger} \mathbf{g}_{\parallel},
\end{equation*}
where $ {\dagger} $ denotes the pseudo-inverse. If 
$ \| \mathbf{v}_{\parallel}  (0) \|_2 \le \delta $, then 
\begin{equation*}
			\label{eq:v_par_solsing}
			\mathbf{v}^*_{\parallel} = \mathbf{v}_{\parallel} \left( 0 \right) = - \Lambda^{\dagger}\mathbf{g}_{\parallel}
\end{equation*}
satisfies all optimality conditions (with $ \sigma^*_{\parallel} = 0
$). Otherwise, i.e., if either $[\mathbf{g}_{\parallel}]_i \ne 0 $ for some
$ 1 \le i \le r $ or $\| \Lambda^{\dagger} \mathbf{g}_{\parallel} \|_2 >
\delta $, then $ \mathbf{v}^*_{\parallel} $ is computed using
\eqref{eq:v_secularsoln}, where $ \sigma^*_{\parallel} $ solves the secular
equation in \eqref{eq:secular} (discussed below).
		
\bigskip

\noindent \textbf{Case 3 ($ \lambda_1 < 0 $):} If $ \mathbf{g}_{\parallel} $ is in the range of $\Lambda - \lambda_1 \mathbf{I}$, 
i.e., $[\mathbf{g}_{\parallel}]_i = 0$ for $ 1 \le i \le r$, then 
we set $ \sigma_{\parallel} = -\lambda_1 $ and 
\begin{equation*}
	\mathbf{v}_{\parallel}\left( -\lambda_1 \right) = -\left( \Lambda -\lambda_1\mathbf{I} \right)^{\dagger} \mathbf{g}_{\parallel}.
\end{equation*}
If 
$\| \mathbf{v}_{\parallel} ( -\lambda_1 )\|_2 \le \delta $, then the solution is given by
\begin{equation}
	\label{eq:v_par_solnhc}
	\mathbf{v}^*_{\parallel} = \mathbf{v}_{\parallel} \left( -\lambda_1 \right) + \alpha  \mathbf{e}_1,
\end{equation} 
where $ \alpha = \sqrt{ \delta^2 -  \left\| \mathbf{v}_{\parallel}\left( -\lambda_1 \right)  \right\|^2_2  } $. 
(This case is referred to as the ``hard case''~\cite{ConGT00a,MorS83}.)
Note that $\mathbf{v}_{\parallel}^*$ satisfies the first optimality condition \eqref{eq:vpar_c1}:
\begin{align*}
	\left( \Lambda  -\lambda_1  \mathbf{I} \right) \mathbf{v}^*_{\parallel} 
	&= \left( \Lambda  -\lambda_1  \mathbf{I} \right) \left( \mathbf{v}_{\parallel} \left( -\lambda_1 \right) + \alpha  \mathbf{e}_1 \right) 
	= -\mathbf{g}_{\parallel}.
\end{align*} 
The second optimality condition \eqref{eq:vpar_c2} is satisfied by observing that
\begin{equation*}
		\| \mathbf{v}^*_{\parallel} \|^2_2 =  \| \mathbf{v}_{\parallel} ( -\lambda_1 )  \|^2_2 + \alpha^2 = \delta^2.
\end{equation*}
Finally, since $ \sigma^*_{\parallel} = -\lambda_1 > 0 $ the other optimality
conditions are also satisfied.

On the other hand, if $ [ \mathbf{g}_{\parallel}
]_i \ne 0$ for some $ 1 \le i \le r $ or $ \| ( \Lambda
-\lambda_1\mathbf{I} )^{\dagger}\mathbf{g}_{\parallel} \|_2 > \delta $,
then $ \mathbf{v}^*_{\parallel} $ is computed using
\eqref{eq:v_secularsoln}, where $ \sigma^*_{\parallel} $ solves the secular
equation \eqref{eq:secular}.

\medskip

\bigskip

\noindent \textbf{The secular equation.}  We now summarize how to find a
solution of the so-called \emph{secular equation}.  Note that from
(\ref{eq:v_sigma}),
\begin{equation*}
	\| \mathbf{v}_{\parallel} (\sigma_{\parallel}) \|_2^2 =  \sum_{i=1}^{m} \frac{(\mathbf{g}_{\parallel})_i^2}{(\lambda_i + \sigma_{\parallel})^2}.
\end{equation*}
If we combine the terms above that correspond to the same eigenvalues and remove the terms with zero 
numerators, then for $\sigma_{\parallel} \ne -\lambda_i$, we have 
$$
	\| \mathbf{v}_{\parallel} (\sigma_{\parallel}) \|_2^2 =  \sum_{i=1}^{\ell} \frac{\bar{a}_i^2}{(\bar{\lambda}_i + \sigma_{\parallel})^2},
$$
where $\bar{a}_i \ne 0$ for $i = 1, \dots, \ell$ and $\bar{\lambda}_i$ are \emph{distinct} eigenvalues of $\mathbf{B}$ with 
$\bar{\lambda}_1 < \bar{\lambda}_2 < \cdots < \bar{\lambda}_{\ell}$.  
Next, we define the function
\begin{equation}
	\phi_{\parallel}\left( \sigma_{\parallel} \right)
	=
	\begin{cases}
		 \displaystyle \frac{1}{  
		 \sqrt{\displaystyle \sum_{i=1}^{\ell} \frac{\bar{a}_i^2}{(\bar{\lambda}_i + \sigma_{\parallel})^2}}
		 } - \frac{1}{\delta} 
		 & \text{if $\sigma_{\parallel} \ne - \bar{\lambda}_i$ where $1 \le i \le \ell$} \\[1.3cm]
		 \displaystyle -\frac{1}{\delta} 
		 & \text{otherwise}.
	\end{cases}
\end{equation}
From the optimality conditions \eqref{eq:vpar_c2}
and \eqref{eq:vpar_c5}, if $\sigma^*_{\parallel} \ne 0$, then
$\sigma_{\parallel}^*$ solves the \emph{secular equation}
\begin{equation}
	\label{eq:secular}
	\phi_{\parallel}\left( \sigma_{\parallel} \right) 
	= 
	0, 
\end{equation}
with $\sigma_{\parallel} \ge \max \{0, -\lambda_1\}$.
Note that $\phi_{\parallel}$ is monotonically increasing 
and concave on the interval $[
- \lambda_1, \infty)$; thus, with a judicious choice of initial $\sigma_{\parallel}^0$, Newton's method can be used to efficiently
compute $\sigma_{\parallel}^*$ in \eqref{eq:secular} (see~\cite{BruEM15}).

\bigskip

	More details on the solution method for subproblem \eqref{eq:sp_vpar_2}
are given in \cite{BruEM15}.

\subsection{Computing $\mathbf{p}^*$}
Given $\mathbf{v}^* = [ \mathbf{v}^*_{\parallel} \,\,
			\mathbf{v}^*_{\perp}]^T$,
the solution to the trust-region subproblem \eqref{eq:trustProblem} 
using either the $(\mathbf{P},2)$ or the $(\mathbf{P},\infty)$ norms
is
\begin{equation}
	\label{eq:p_general}
		\mathbf{p}^* = \mathbf{P}\mathbf{v}^* = \mathbf{P}_{\parallel}\mathbf{v}^*_{\parallel}+\mathbf{P}_{\perp}\mathbf{v}^*_{\perp}.
\end{equation}
(Recall that using either of the two norms generates the same $ \bfm{v}_\perp^*$
but different $ \bfm{v}_\parallel^*$.)  It remains to show how to form
$\mathbf{p}^*$ in (\ref{eq:p_general}).  Matrix-vector products
  involving $ \mathbf{P}_{\parallel} $ are possible using \eqref{eqn-Pparallel} or \eqref{eqn-Pparallel_nc}, and thus, $ 
  \mathbf{P}_{\parallel}\mathbf{v}^*_{\parallel} $ can be computed;
  however, an explicit formula to compute products
 $\mathbf{P}_{\perp}$ is
  not available.  To compute the second term,
$\mathbf{P}_{\perp}\mathbf{v}^*_{\perp}$, we observe that $
\mathbf{v}^*_{\perp} $, as given in \eqref{eq:soln_vperp}, is a multiple of
either $ \mathbf{g}_{\perp} = \mathbf{P}^T_{\perp}\mathbf{g}$ or a vector $
\mathbf{u}$ with unit length, depending on the sign of $\gamma$ and
the magnitude of $\mathbf{g}_\perp$.  In the latter case, define $ \mathbf{u} = \frac{
  \mathbf{P}^T_{\perp}\mathbf{e}_i }{ \| \mathbf{P}^T_{\perp}\mathbf{e}_i
  \|_2 } $, where $ i \in \left\{1,2,\ldots,k+2 \right\} $ is the first
index such that $ \left\| \mathbf{P}^T_{\perp}\mathbf{e}_i \right\|_2 \neq
0 $.  (Such an $\mathbf{e}_i$ exists since rank($\mathbf{P}_{\perp}) =
n-m$.)   Thus, we obtain
\begin{equation}
	\label{eq:opt_p}
	\mathbf{p}^* = \mathbf{P}_{\parallel}(\mathbf{v}^*_{\parallel} - \mathbf{P}^T_{\parallel}\mathbf{w}^*) + \mathbf{w}^*,
\end{equation}
where 
\begin{equation}
	\label{eq:opt_w}
	\mathbf{w}^* = 
	\begin{cases}			
	-\frac{1}{\gamma}\mathbf{g} 														& \text{ if } \gamma > 0 \text{ and } \left \| \mathbf{g}_{\perp} \right \|_2 \le \delta
|\gamma |, \\
	\frac{\delta}{\left\| \mathbf{P}^T_{\perp}\mathbf{e}_i \right\|_2}\mathbf{e}_i		& \text{ if } \gamma \le 0 \text{ and } \| \mathbf{g}_{\perp} \|_2 = 0, \\
	-\frac{ \delta }{\| \mathbf{g}_{\perp} \|_2}\mathbf{g}   					& \text{ otherwise.}			
	\end{cases}
\end{equation}

Algorithm~\ref{alg-w} summarizes the computation of $\mathbf{w}^*$.
 \begin{algorithm}[h!]
   \caption{Computing $\mathbf{w}^*$}\label{alg-w}
 \begin{algorithmic}[1]
\ENSURE [$\mathbf{w}^*$, $\beta$, $\text{hasBeta}$]=\texttt{ComputeW}($\mathbf{g},\delta,\gamma,\|\mathbf{g}_\perp\|_2, \mathbf{\Pi},\mathbf{\Psi},\mathbf{R}_\ddagger,\mathbf{U},\tau,[\jb{\texttt{varargin}=\{\bfm{S}, \bfm{Y} \}}]$);
\IF {$\gamma>0$ \AND $\|\mathbf{g}_\perp\|_2\le \delta\gamma$} 
\STATE {$\beta \gets -(1/\gamma)$, $ \text{hasBeta} \gets 1 $;}
\STATE {$\mathbf{w}^*\gets \beta\mathbf{g}$;}
\ELSIF{$\gamma\le 0$ \AND $\|\mathbf{g}_\perp\|_2<\tau$} \STATE
      {Find the first index $i$ such that $\|\mathbf{P}^T_\perp \mathbf{e}_i\|_2\ne 0$;}
      \STATE{$\beta \gets 0$, $ \text{hasBeta} \gets 0 $;}
      \STATE {$\mathbf{w}^*\gets \left(\delta/\|\mathbf{P}^T_\perp \mathbf{e}_i\|_2\right)\mathbf{e}_i$;}
\ELSE
\STATE {$\beta \gets -(\delta/\|\mathbf{g}_\perp\|_2)$, $ \text{hasBeta} \gets 1 $;}
\STATE {$\mathbf{w}^*\gets \beta\mathbf{g}$;}
\ENDIF
\RETURN $\mathbf{w}^*$;
 \end{algorithmic}
 \end{algorithm}

The quantities $ \left\| \mathbf{g}_{\perp} \right\|_2 $ and $\left\| \mathbf{P}^T_{\perp}\mathbf{e}_i \right\|_2$
are computed using the orthogonality of $\bfm{P}$, which implies
\begin{equation}
	\label{eq:ng_perp}
	\left \| \mathbf{g}_{\parallel} \right \|^2_2 +  \| \mathbf{g}_{\perp}  \|^2_2 =  \| \mathbf{g}  \|^2_2, 
	\words{and} 
	\| \mathbf{P}^T_{\parallel}\mathbf{e}_i  \|^2_2 + \| \mathbf{P}^T_{\perp}\mathbf{e}_i \|^2_2 = 1.
\end{equation}
Then $\| \mathbf{g}_{\perp} \|_2 = \sqrt{ \| \mathbf{g}\|^2_2 - \|
  \mathbf{g}_{\parallel} \|^2_2 } $ and $ \|
\mathbf{P}^T_{\perp}\mathbf{e}_i \|_2 = \sqrt{ 1 -
  \|\mathbf{P}^T_{\parallel}\mathbf{e}_i \|^2_2 } $.  Note that
$\mathbf{v}_{\perp}^*$ is never explicitly computed. 
\jb{Since $ \bs{\Psi} $ is either explicitly computed when a constant initialization is used or represented through
$ \bfm{S} $ and $ \bfm{Y} $, the optional input $ [\texttt{varargin}= \{\bfm{S}, \bfm{Y} \}]$ can be used to
pass $ \bfm{S}, \bfm{Y} $ if $ \bs{\Psi} $ is represented implicitly.}

\section{The Proposed Algorithms}
In this section, we summarize Section 3 in two algorithms that solve 
the trust-region subproblem using the $(\mathbf{P},\infty)$
and the $(\mathbf{P},2)$ norms.  The required inputs \jb{depend on the initialization strategy and often} include
$\mathbf{g}$, $\mathbf{S}$,
$\mathbf{Y}$, $\gamma$, and $\delta$ which define the trust-region subproblem (including the \LSR{} matrix).  The input $\tau$ is a small positive number
used as a tolerance.  The output of each algorithm is $\bfm{p}^*$, the solution
to the trust-region subproblem in the given shape-changing norm.
In Algorithm~\ref{alg-pinfty}, we detail the algorithm for solving
(\ref{eq:trustProblem}) using the $(\mathbf{P},\infty)$ norm to define
the subproblem; Algorithm~\ref{alg-p2} solves the subproblem using the
$(\mathbf{P},2)$ norm.

Both algorithms accept either the matrices that hold the quasi-Newton
pairs,
$\mathbf{S}$ and $\mathbf{Y}$, or factors for the compact formulation $\mathbf{\Psi}$ and
$\mathbf{M}^{-1}$.  To reflect this option, the second and third input
parameters are labeled ``$S\lor\Psi$'' and ``$Y\lor M^{-1}$'' in both algorithms.
If the inputs are the quasi-Newton pairs $ \bfm{S}, \bfm{Y} $, the input parameter
\texttt{flag} must be ``0'', and then factors for the compact
formulation are computed; if the inputs are factors for the compact
formulation,
\texttt{flag} must be ``1'', and then $\mathbf{\Psi}$ and $\mathbf{M}^{-1}$ are
set to the second and third inputs.
\jb{Another option is to pass the product $ \bs{\Psi}^T \bs{\Psi} $ and the matrix $ \bfm{M}^{-1} $ along with
$ \bfm{S} $ and $ \bfm{Y} $. This can be particularly advantageous when the matrix $ \bs{\Psi}^T \bs{\Psi} $ is
updated, instead of recomputed (by using e.g., Procedure 2).}

\begin{algorithm}[H]
\caption{The SC-SR1 algorithm for the 
 shape-changing $(\mathbf{P},\infty)$ norm}
\label{alg-pinfty}
\begin{algorithmic}[1]
\ENSURE [$\mathbf{p}^*$]=\texttt{sc\_sr1\_infty}($\mathbf{g}$, $S\lor\Psi$,
$Y\lor M^{-1}$, $\gamma$, $\delta$, flag, $\jb{[\texttt{varargin}=\{\Psi^T\Psi, M^{-1} \}]}$)
\STATE{Pick $\tau$ such that $0<\tau\ll 1$;}
\IF{flag$=0$}  
	\STATE{\jb{$\mathbf{S}\gets S\lor\Psi$ and $\mathbf{Y}\gets
		  Y\lor M^{-1}$;}}
	\IF{\jb{$ \texttt{isempty(varargin)}$}}
		\STATE {\jb{Compute $\mathbf{\Psi}$ and $\mathbf{M}^{-1}$ as in
		  (\ref{eqn-PsiM})}}
	\ELSE
		\STATE{ \jb{$\Psi \lor \Psi^T \Psi \gets \texttt{varargin}\{1\}$ and
		$\bfm{M}^{-1} \gets \texttt{varargin}\{2\}$;}}
	\ENDIF
\ELSE
	 \STATE {$\mathbf{\Psi}\gets S\lor\Psi$; $\mathbf{M}^{-1}\gets Y\lor M^{-1}$;}
\ENDIF
\STATE  [$\mathbf{R}_\ddagger$,
$\mathbf{\Lambda}$, $\mathbf{U}$,
  $\mathbf{\Pi}$, $J$]=\texttt{ComputeSpectral}($\Psi \lor \Psi^T \Psi$, $\mathbf{M}^{-1}$, $\gamma$, $\tau$);
\STATE{\jeo{$m\gets|J|$ };}
\STATE $\mathbf{g}_\parallel\gets \bfm{P}_\parallel^T \bfm{g}$ using (\ref{eqn-Pparallel});
\STATE{$\|\mathbf{g}_\perp\|\gets\sqrt{\|\mathbf{g}\|_2^2-\|\mathbf{g}_\parallel\|^2}$;}
\IF{$\|\mathbf{g}_\perp\|<\tau$} 
\STATE{$\|\mathbf{g}_\perp\|\gets 0$;}
\ENDIF
\FOR{$i=1$ \TO $m$}\STATE{\IF{$\left|[\mathbf{g}_\parallel]_i\right|<\delta\left|[\mathbf{\Lambda}]_{ii}\right|$ \AND
$[\mathbf{\Lambda}]_{ii}>\tau$} 
\STATE{$[\mathbf{v}_\parallel]_i\gets -[\mathbf{g}_\parallel]_i/[\mathbf{\Lambda}]_{ii}$;}
\ELSIF{$\left|[\mathbf{g}_\parallel]_i\right|<\tau$ \AND $\left|[\mathbf{\Lambda}]_{ii}\right|<\tau$} 
\STATE{$[\mathbf{v}_\parallel]_i\gets\delta/2$;}
\ELSIF{$\left|[\mathbf{g}_\parallel]_i\right|>\tau$ \AND $\left|[\mathbf{\Lambda}]_{ii}\right|<\tau$} 
\STATE{$[\mathbf{v}_\parallel]_i\gets -\text{sgn}\left([\mathbf{g}_\parallel]_i\right)\delta$;}
\ELSIF{$\left|[\mathbf{g}_\parallel]_i\right|<\tau$ \AND $[\mathbf{\Lambda}]_{ii}<-\tau$} 
\STATE{$[\mathbf{v}_\parallel]_i\gets\delta$;}
\ELSE
\STATE{$[\mathbf{v}_\parallel]_i\gets-\left(\delta/\left|[\mathbf{g}_\parallel]_i\right|\right) [\mathbf{g}_\parallel]_i$;}
\ENDIF}
\ENDFOR
\STATE{[$\mathbf{w}^*$,$\beta$,$\text{hasBeta}$]=\texttt{ComputeW}($\mathbf{g},\delta,\gamma,\|\mathbf{g}_\perp\|_2, \mathbf{\Pi},\mathbf{\Psi},\mathbf{R}_\ddagger,\mathbf{U},\tau,\jb{[\texttt{varargin}=\{\bfm{S},\bfm{Y}\}]} $);  }
\IF{\text{hasBeta} = 1}
\STATE{$\mathbf{p}^*\gets \mathbf{P}_\parallel ( \bfm{v}_\parallel- \beta \mathbf{g}_\parallel) + \mathbf{w}^*$}
\ELSE
\STATE{$\mathbf{p}^*\gets \mathbf{P}_\parallel ( \bfm{v}_\parallel- \mathbf{P}_\parallel^T\mathbf{w}^*) + \mathbf{w}^*$}
\ENDIF
\RETURN $\mathbf{p}^*$
\end{algorithmic}
\end{algorithm}

The computation of $\mathbf{p}^*$ in both Algorithms~\ref{alg-pinfty}
and~\ref{alg-p2} is performed as in (\ref{eq:opt_p}) using two
matrix-vector
products with $\mathbf{P}^T_\parallel$ and $\mathbf{P}_\parallel$,
respectively, in order to avoid matrix-matrix products.
Products with $\mathbf{P}_\parallel$ are done using the factors
$\mathbf{\Psi}$, $\mathbf{R}$,
        and $\mathbf{U}$ (see Section~\ref{sec-eigs}).

%

\begin{algorithm}[H]
\footnotesize
\caption{The SC-SR1 algorithm for the 
 shape-changing $(\mathbf{P},2)$ norm}
\label{alg-p2}
\begin{algorithmic}[1] 
\ENSURE [$\mathbf{p}^*$]=\texttt{sc\_sr1\_2}($\mathbf{g}$, $S\lor\Psi$, $Y\lor M^{-1}$, $\gamma$, $\delta$, flag, $\jb{[\texttt{varargin}=\{\Psi^T\Psi, M^{-1} \}]}$)
\STATE{Pick $\tau$ such that $0<\tau\ll 1$;}
\IF{flag$=0$}  
	\STATE{\jb{$\mathbf{S}\gets S\lor\Psi$ and $\mathbf{Y}\gets
		  Y\lor M^{-1}$;}}
	\IF{\jb{$ \texttt{isempty(varargin)}$}}
		\STATE {\jb{Compute $\mathbf{\Psi}$ and $\mathbf{M}^{-1}$ as in
		  (\ref{eqn-PsiM})}}
	\ELSE
		\STATE{ \jb{$\Psi \lor \Psi^T \Psi \gets \texttt{varargin}\{1\}$ and
		$\bfm{M}^{-1} \gets \texttt{varargin}\{2\}$;}}
	\ENDIF
\ELSE
	 \STATE {$\mathbf{\Psi}\gets S\lor\Psi$; $\mathbf{M}^{-1}\gets Y\lor M^{-1}$;}
\ENDIF
  \STATE  [$\mathbf{R}_\ddagger$,
$\mathbf{\Lambda}$, $\mathbf{U}$,
$\mathbf{\Pi}$, $J$]=\texttt{ComputeSpectral}($\Psi \lor \Psi^T \Psi$, $\mathbf{M}^{-1}$, $\gamma$, $\tau$);
\STATE{$\mathbf{g}_\parallel\gets \bfm{P}_\parallel^T \bfm{g}$ using (\ref{eqn-Pparallel}), and $\|\mathbf{g}_\perp\|\gets\sqrt{\|\mathbf{g}\|_2^2-\|\mathbf{g}_\parallel\|^2}$ };
\IF{$\|\mathbf{g}_\perp\|<\tau$}
\STATE{$\|\mathbf{g}_\perp\|\gets 0$;}
\ENDIF
\IF{$[\mathbf{\Lambda}]_{11}>\tau$}
\IF{$\|\mathbf{\Lambda^{-1}}\mathbf{g}_\parallel\|<\delta$}
\STATE{$\sigma_\parallel\gets 0$;}
\ELSE
\STATE{Use Newton's method with $\sigma_0=0$ to find $\sigma_\parallel$, a solution to (\ref{eq:v_sigma});}
\ENDIF

\STATE{$\mathbf{v}_\parallel \gets -\left(\mathbf{\Sigma}+\sigma_\parallel \bfm{I}\right)^{-1}\mathbf{g}_\parallel$;}
\ELSIF{$\left|[\mathbf{\Lambda}]_{11}\right|<\tau$ }
  \STATE{Define $r$ to be the first $i$ such that $\left|\Lambda_{ii}\right|>\tau$;}
   \IF{$\left| \mathbf{g}_{ii}\right|<\tau $ for $1\le i \le r$ \AND $\|\mathbf{\Lambda}^\dagger \mathbf{g}_\parallel\|<\delta$}
       \STATE{$\sigma_\parallel\gets 0$;}
       \STATE{$\mathbf{v}_\parallel\gets -\mathbf{\Lambda}^\dagger \mathbf{g}_\parallel$;}
       \ELSE
       \STATE {$\hat{\sigma}=\text{max}_i ( [\mathbf{g}_\parallel]_i/\delta-\mathbf{\Lambda}_{ii})$;}   
      \STATE{Use Newton's method with $\sigma_0=\hat{\sigma}$ to find $\sigma_\parallel$, a solution to (\ref{eq:v_sigma});}
     \STATE{$\mathbf{v}_\parallel \gets -\left(\mathbf{\Lambda}+\sigma_\parallel \bfm{I}\right)^{-1}\mathbf{g}_\parallel$;}
     \ENDIF
     \ELSE  
     \STATE {Define $r$ to be the first $i$ such that $\left|[\Lambda]_{ii}\right|>\tau$;}

     \IF {$\left|\mathbf{g}_{ii}\right|<\tau $ for $1\le i \le r$ }
     \STATE{$\sigma_\parallel=-[\mathbf{\Lambda}]_{11}$, $\mathbf{v} \gets \left(\mathbf{\Lambda}-[\mathbf{\Lambda}]_{11}\bfm{I}\right)^\dagger \mathbf{g}_\parallel$;}

     \IF{$\|\mathbf{v}\|<\delta$}
     \STATE{$\alpha \gets \sqrt{\delta^2-\|\mathbf{v}\|^2}$;}
     \STATE{$\mathbf{v}_\parallel = \bfm{v}+\alpha \bfm{e}_1$, where $\bfm{e}_1$ is the
       first standard basis vector;}
     \ELSE
     \STATE{$\hat{\mathbf{\sigma}}\gets   \text{max}_i (
       [\mathbf{g}_\parallel]_i/\delta-[\mathbf{\Lambda}]_{ii})$;}
     \STATE{Use Newton's method with $\sigma_0=\text{max}(\hat{\sigma},0)$ to find
       $\sigma_\parallel$, a solution to (\ref{eq:v_sigma});}
     \STATE{$\mathbf{v}_\parallel \gets       -\left(\mathbf{\Lambda}+\sigma_\parallel \bfm{I}\right)^{-1}\mathbf{g}_\parallel$;}
     \ENDIF
     \ELSE
    \STATE{$\hat{\mathbf{\sigma}}\gets   \text{max}_i (
      [\mathbf{g}_\parallel]_i/\delta-[\mathbf{\Lambda}]_{ii})$;}
    \STATE{Use Newton's method with $\sigma_0=\text{max}(\hat{\sigma},0)$ to find
      $\sigma_\parallel$, a solution to (\ref{eq:v_sigma});}
    \STATE{$\mathbf{v}_\parallel \gets
      -\left(\mathbf{\Lambda}+\sigma_\parallel
        \bfm{I}\right)^{-1}\mathbf{g}_\parallel$;}
\ENDIF
\ENDIF
\STATE{[$\mathbf{w}^*$,$\beta$,$\text{hasBeta}$]=\texttt{ComputeW}($\mathbf{g},\delta,\gamma,\|\mathbf{g}_\perp\|_2, \mathbf{\Pi},\mathbf{\Psi},\mathbf{R}_\ddagger,\mathbf{U},\tau,\jb{[\texttt{varargin}=\{\bfm{S},\bfm{Y}\}]} $); }
\IF{\text{hasBeta} = 1}
\STATE{$\mathbf{p}^*\gets \mathbf{P}_\parallel ( \bfm{v}_\parallel- \beta \mathbf{g}_\parallel) + \mathbf{w}^*$}
\ELSE
\STATE{$\mathbf{p}^*\gets \mathbf{P}_\parallel ( \bfm{v}_\parallel- \mathbf{P}_\parallel^T\mathbf{w}^*) + \mathbf{w}^*$}
\ENDIF
\end{algorithmic}
\end{algorithm}

\jb{Besides the optional arguments,} the {\small MATLAB} implementation of both algorithms have an
additional input and output variable.  The additional input variable
is a verbosity setting; the additional output variable is a flag that
reveals whether there were any detected run-time errors.

\jb{As described in Section \ref{sec-lm} the limited-memory updating techniques vary with the choice of 
initialization strategy $ \bfm{B}^{(k)}_0 = \gamma_k \bfm{I} $. A commonly used value is $ \gamma_k = \frac{ \| \bk{y} \|^2 }{ \bk{s}^T \bk{y} } $ \cite{bb88}. Recall from Section \ref{sec-eigs} that $ \gamma_k $ is the eigenvalue for the large $ n-m $ dimensional subspace,
spanned by the eigenvector in $ \bfm{P}_{\perp} $. At a local minimum all eigenvalues of $ \nabla^2 f(\bfm{x}) $ will be non-negative,
motivating non-negative values of $\gamma_k$. For our implementation we tested three different strategies, one of which
uses a constant initialization (C Init.)
\begin{equation}
	\label{eq:init}
	\gamma_k =
	\begin{cases}
		\text{max}\left(\text{min}\left(\frac{\| \bfm{y}_0 \|^2}{\bfm{s}^T_0 \bfm{y}_0}, \gamma_{\text{max}}\right),1\right) 	& \text{C Init.} \\
		\frac{\| \bk{y} \|^2}{\bk{s}^T \bk{y}} \: \left(\text{if } \bk{s}^T \bk{y} > 0 \right)					& \text{Init. 1} \\ 
		\text{max}\left(\frac{\| \bfm{y}_{k-q} \|^2}{\bfm{s}^T_{k-q} \bfm{y}_{k-q}}, \cdots, \frac{\| \bk{y} \|^2}{\bk{y}^T \bk{s}} \right) & \text{Init. 2}
	\end{cases}
\end{equation} 
Observe that Init. 2 includes the additional parameter $ q > 0 $, which determines the number of pairs $\{\bfm{s}_i,\bfm{y}_i\}$ to use.
For C Init., the parameter $\gamma_{\text{max}}$ ensures that the constant initialization, which uses $ \bfm{s}_0, \bfm{y}_0 $ does not
exceed this threshold value. In the experiments the parameter is set as  $ \gamma_{\text{max}} = 1 \times 10^4 $.}

\subsection{Computational Complexity}
\label{subsec:cc}

We estimate the cost of one iteration using the proposed method to solve
the trust-region subproblem defined by shape-changing norms (\ref{eq:sc_2})
and (\ref{eq:sc_inf}).  We make the practical assumption that
$\gamma>0$.
Computational savings can be achieved by reusing 
previously computed matrices and not forming certain matrices explicitly.
We begin by highlighting the \jb{case when a non-constant initialization strategy is used.} 
First, we do not form $\mathbf{\Psi}=\mathbf{Y}-\gamma\mathbf{S}$ explicitly.  Rather, we
compute matrix-vector products with
$\mathbf{\Psi}$ by computing matrix-vector
products with $\mathbf{Y}$ and $\mathbf{S}$.  
Second, to form $\mathbf{\Psi}^T\mathbf{\Psi}$, 
we only store and update the small $m \times m$ matrices
$\mathbf{Y}^T\bfm{Y}$, $\mathbf{S}^T\bfm{Y}$, and $\mathbf{S}^T\bfm{S}$. 
This update involves only $3m$ vector inner products.
Third, assuming we have already obtained the Cholesky
factorization of $\mathbf{\Psi}^T\mathbf{\Psi}$ associated with the
previously-stored limited-memory pairs, it is possible to update the
Cholesky factorization of the new $\mathbf{\Psi}^T\mathbf{\Psi}$ at a cost
of $O(m^2)$ \cite{Ben65,GilGMS74}.

We now consider the dominant cost for a single subproblem solve.
The eigendecomposition $\mathbf{R_{\dagger}\Pi}^T \bfm{M\Pi R_{\dagger}} = \mathbf{U}
\mathbf{\hat{\Lambda}} \mathbf{U}^T$ costs $O(m^3) = \left ( \frac{m^2}{n}
\right ) O(mn)$, where $m \ll n$.  To compute $\mathbf{p}^*$ in
\eqref{eq:opt_p}, one needs to compute $\mathbf{v}^*$ from Section
\ref{sec:v_parallel} and $\mathbf{w}^*$ from \eqref{eq:opt_w}.
The dominant cost for computing $\mathbf{v}^*$ and
$\mathbf{w}^*$ is forming 
$\mathbf{\Psi}^T\mathbf{g}$, which 
requires $\jb{2mn}$ operations.  
\jb{Note that both $ \bfm{v}^*_{\parallel} $ and $ \bfm{P}^T_{\parallel} \bfm{w}^* $ are typically computed from $ \bfm{P}^T_{\parallel} \bfm{g} $,
whose main operation is $ \bs{\Psi}^T \bfm{g} $. Subsequently, computing $ \bfm{P}_{\parallel}( \bfm{v}^*_{\parallel} - \bfm{P}_{\parallel}^T \bfm{w}^*) $ incurs $O(2mn)$ additional multiplications, as this operation reduces to $ \bs{\Psi} \bfm{f} $ for a vector $ \bfm{f} $. Thus,
the dominant complexity is $ O(2mn + 2mn) = O(4mn) $.}
The following theorem summarizes the dominant computational costs.

\begin{theorem} 
  The dominant computational cost of solving one trust-region subproblem for the
  proposed method is $4mn$ floating point operations.
\end{theorem}

We note that the floating point operation count of $O(4mn)$  
is the same cost as for \LBFGS~\cite{Noc80}. 

\jb{If a constant initialization is used the complexity can essentially be halved, because the mat-vec applies $ \bs{\Psi}^T \bfm{g} $
and $ \bs{\Psi} \bfm{f} $ (for some vector $ \bfm{f} $) each take $ O(mn) $ multiplications for a total of $ O(2mn) $.}

\subsection{Characterization of global solutions}

It is possible to characterize global solutions to the trust-region subproblem defined by shape-changing norm
$\left(\mathbf{P},2\right)$-norm.  The following theorem is based on well-known
optimality conditions for the 
two-norm trust-region subproblem~\cite{Gay81,MorS83}.  
\begin{theorem} \label{thrm-charact} A vector 
 $\mathbf{p}^*\in \mathbb{R}^n$
such that $		\left\| \mathbf{P}^T_{\parallel} \mathbf{p}^* \right\|_2 \le \delta$ and  
$	\left\| \mathbf{P}^T_{\perp} \mathbf{p}^* \right\|_2\le \delta,$
 is a global solution of 
(\ref{eq:trustProblem}) defined by the $\left(\mathbf{P},2 \right)$-norm if
and only if there exists unique $\sigma_\parallel^*\ge 0$ and
$\sigma_\perp^*\ge 0$ such that
\begin{eqnarray}		\label{eq:optim_sc2} 
\left( \mathbf{B} + \mathbf{C}_{\parallel} \right)\mathbf{p}^* +\mathbf{g} = 0,
\quad \sigma^*_{\parallel}  \left( \left\| \mathbf{P}^T_{\parallel}\mathbf{p}^* \right \|_2 - \delta \right) = 0, \quad	
		\sigma^*_{\perp}  \left( \left \|\mathbf{P}_\perp^T\mathbf{p}^*\right\|_2 - \delta \right) = 0, 
\begin{array}{lcrrcl}		
\end{array}
\end{eqnarray}
where $ \mathbf{C}_{\parallel} \defined \sigma^*_{\perp}\mathbf{I} 
+ \left( \sigma^*_{\parallel} - \sigma^*_{\perp} \right)\mathbf{P}_{\parallel}\mathbf{P}^T_{\parallel}$,
the matrix $\mathbf{B+C_\parallel}$ is positive semi-definite,
and $ \mathbf{P}=[\mathbf{P}_\parallel \,\, \mathbf{P}_\perp]$
and
$\Lambda=\diag(\lambda_1,\ldots,\lambda_{m})=\hat{\Lambda}+\gamma\mathbf{I}$
are as in (\ref{eqn-PL}).
\end{theorem}	

When run in the ``verbose'' mode,
\texttt{sc\_sr1\_2.m} returns values needed to
establish the optimality of $\mathbf{p}^*$
using this theorem.  In particular, the code computes 
$\sigma^*_\parallel$, which, depending on the case, is either 0, the
absolute value of the most negative
eigenvalue, or obtained from Newton's method.  The code
also computes $\sigma^*_\perp$ using
(\ref{eq:soln_sigmaperp}), and $\|\mathbf{P}^T_\perp \bfm{p}^*\|$
is computed by noting that  $\|\mathbf{P}^T_\perp \bfm{p}^*\|_2^2=
\|\mathbf{p}^*\|^2_2-\|\mathbf{P}^T_\parallel \bfm{p}^*\|^2_2$.
The variables \texttt{opt1}, \texttt{opt2}, and \texttt{opt3}
contain the errors in each of the equations 
in (\ref{eq:optim_sc2}); \texttt{spd\_check} finds the minimum 
eigenvalue of $(\mathbf{B}+\mathbf{C}_\parallel)$ in
(\ref{eq:optim_sc2}), 
enabling one to ensure $(\mathbf{B}+\mathbf{C}_\parallel)$ is positive
definite; and $\sigma^*_\parallel$ and
$\sigma^*_\perp$ are displayed to verify that they are nonnegative.



\section{Numerical experiments}
In this section, we report on numerical experiments with the proposed
shape-changing SR1 (\SCSR) algorithm implemented in \MATLAB{} to solve
limited-memory SR1 trust-region subproblems. 
\jb{The experiments are divided into solving the TR subproblems
with Algorithms \ref{alg-pinfty} and \ref{alg-p2}, and general unconstrained minimization problems, which use
the TR subproblem solvers, using 62 large-scale CUTEst problems \cite{GouOT03}.}

\subsection{$(\mathbf{P},2)$-norm results}
The \SCSR{} algorithm
  was tested on randomly-generated problems of size $n=10^3$ to $n=10^7$,
  organized as five experiments when there is no closed-form solution to
  the shape-changing trust-region subproblem and one experiment designed to
  test the \SCSR{} method in the so-called ``hard case''.  These six
  cases only occur using the
  $(\mathbf{P},2)$-norm trust region.  (In the case of the
  $(\mathbf{P},\infty)$ norm, $\mathbf{v}_\parallel^*$ has the closed-form
  solution given by (\ref{eq:vpar_inf}).) 
The six experiments are outlined as follows:

\begin{enumerate}

\item[(E1)] {$\mathbf{B}$ is positive definite with $\| \mathbf{v}_\parallel(0)
  \|_2 \ge \delta$}.
\item[(E2)] {$\mathbf{B}$ is positive semidefinite and singular with
$[\mathbf{g}_{\parallel}]_i \ne 0 $ for some $ 1 \le i \le r $.}
\item[(E3)] {$\mathbf{B}$ is positive semidefinite and singular with 
$[\mathbf{g}_{\parallel}]_i = 0 $ for $ 1 \le i \le r $ and 
$\| \Lambda^{\dagger} \mathbf{g}_{\parallel} \|_2 > \delta $.}
\item[(E4)] {$\mathbf{B}$ is indefinite and 
$[\mathbf{g}_{\parallel}]_i = 0 $ for $ 1 \le i \le r $ with
$ \| ( \Lambda
-\lambda_1\mathbf{I} )^{\dagger}\mathbf{g}_{\parallel} \|_2 > \delta $.}
\item[(E5)] {$\mathbf{B}$ is indefinite and 
$ [ \mathbf{g}_{\parallel}
]_i \ne 0$ for some $ 1 \le i \le r $ .}
\item [(E6)] {$\mathbf{B}$ is indefinite and 
$[\mathbf{g}_{\parallel}]_i = 0 $ for $ 1 \le i \le r $ with
$\| \mathbf{v}_{\parallel} (
    -\lambda_1 )\|_2 \le \delta $ (the ``hard case'').}
\end{enumerate}

For these experiments, $\mathbf{S}$, $\mathbf{Y}$, and $\mathbf{g}$ were
randomly generated and then altered to satisfy the requirements described
above by each experiment.  In experiments (E2) and (E5), $\delta$ was
chosen as a random number.  (In the other experiments, $\delta$ was
set in accordance with the experiments' rules.)  All randomly-generated vectors and matrices were
formed using the \MATLAB{} \texttt{randn} command, which draws from the
standard normal distribution.  
The initial \SR{} matrix was set to
$\mathbf{B}_0=\gamma \mathbf{I}$, where $\gamma=|10*\texttt{randn(1)}|$.
Finally, the number of limited-memory updates $m$ was set to 5, and $r$
was set to 2.  In the five cases when there is no closed-form solution,
\SCSR{} uses Newton's method to find a root of $\phi_\parallel$. We
  use the same procedure as in \cite[Algorithm 2]{BruEM15} to initialize
  Newton's method since it guarantees monotonic and quadratic convergence
  to $\sigma^*$.  The Newton iteration was terminated when the $i$th
iterate satisfied $\|\phi_\parallel(\sigma^i)\|\le
\texttt{eps}\cdot\|\phi_\parallel(\sigma^0)\| + \sqrt{\texttt{eps}}$, where
$\sigma^0$ denotes the initial iterate for Newton's method and
$\texttt{eps}$ is machine precision.  This stopping criteria is both a
relative and absolute criteria, and it is the only stopping criteria used
by \SCSR.  

\medskip

In order to report on the accuracy of the subproblem solves, we make use of
the optimality conditions found in Theorem~\ref{thrm-charact}.
For each experiment, we report the following:
(i) the norm of the residual of the first optimality condition,
\texttt{opt 1}  $\defined \| (\mathbf{B} + \mathbf{C}_{\parallel}) \mathbf{p}^* +
\mathbf{g} \|_2$; (ii) the first complementarity condition, \texttt{opt
  2}  $\defined |\sigma^*_{\parallel} ( \|\mathbf{P}^T_{\parallel}\mathbf{p}^*
\|_2 - \delta ) |$; (iii) the second complementarity condition,
\texttt{opt 3} $\defined \|\sigma^*_{\perp} ( \|\mathbf{P}^T_{\perp}\mathbf{p}^*
\|_2 - \delta ) |$; (iv) the minimum eigenvalue of
$\mathbf{B}+\mathbf{C}_\parallel$; (v) $\sigma_\parallel^*$;
(vi) $\sigma_\perp^*$; (vii) $\gamma$;
and (viii) time.
The quantities (i)-(vi) are reported to check the
optimality conditions given in Theorem 4.2.
Finally, we ran each experiment five times
and report one representative result for each experiment.

\medskip

\begin{table}[!h]
\setlength\tabcolsep{.95mm}
 \caption{  Experiment 1: $\mathbf{B}$ is positive definite with $\| \mathbf{v}_\parallel(0)
  \|_2 \ge \delta$.} {
  \centering
\begin{tabular}{|c|c|c|c|c|c|c|c|c|}
   \hline $n$   	& \texttt{opt 1} &  \texttt{opt 2} &
\texttt{opt 3} & \texttt{min($\lambda(B+C_\parallel)$)} & 
   $ \sigma^*_{\parallel}$ & $\sigma^*_{\perp} $	& $ \gamma $ &  time  \\ 
\hline
 $1\times 10^3$& \texttt{2.45e-14} & \texttt{0.00e+00} &  \texttt{2.45e-14} & \texttt{4.33e+01} & \texttt{1.09e+01} & \texttt{5.89e+02}& \texttt{1.63e+01} & \texttt{9.97e-03} \\ \hline
  $ 1\times 10^4 $& \texttt{1.21e-13} & \texttt{2.82e-16} &  \texttt{4.26e-13} & \texttt{3.25e+01} & \texttt{8.14e+00} & \texttt{1.98e+03}& \texttt{1.22e+01} & \texttt{1.55e-03} \\ \hline
 $1\times 10^5$& \texttt{5.32e-13} & \texttt{2.28e-16} &  \texttt{1.40e-13} & \texttt{2.19e+01} & \texttt{5.47e+00} & \texttt{5.05e+03}& \texttt{8.14e+00} & \texttt{4.49e-03} \\ \hline
 $1\times 10^6$& \texttt{3.56e-12} & \texttt{5.51e-16} &  \texttt{2.05e-11} & \texttt{1.44e+01} & \texttt{3.61e+00} & \texttt{9.57e+03}& \texttt{5.32e+00} & \texttt{8.03e-02} \\ \hline
 $1\times 10^7$& \texttt{1.46e-11} & \texttt{1.16e-11} &  \texttt{3.64e-11} & \texttt{4.07e+01} & \texttt{1.02e+01} & \texttt{5.52e+04}& \texttt{1.52e+01} & \texttt{9.66e-01} \\ \hline
\end{tabular} }
\end{table}

\setlength\tabcolsep{.95mm}
\begin{table}[!h]
\caption{ Experiment 2: $\mathbf{B}$ is positive semidefinite and singular and
$[\mathbf{g}_{\parallel}]_i \ne 0 $ for some $ 1 \le i \le r $.} {
\centering
\begin{tabular}{|c|c|c|c|c|c|c|c|c|}
   \hline $n$   	& \texttt{opt 1} &  \texttt{opt 2} &
\texttt{opt 3} & \texttt{min($\lambda (B+C_\parallel)$)} & 
$ \sigma^*_{\parallel}$ & $\sigma^*_{\perp} $	& $ \gamma $ &  time  \\
\hline
$1\times 10^3$& \texttt{1.14e-14} & \texttt{0.00e+00} &  \texttt{0.00e+00} & \texttt{9.19e+00} & \texttt{9.19e+00} & \texttt{5.45e+02}& \texttt{1.82e+01} & \texttt{3.24e-03} \\ \hline
  $ 1\times 10^4 $& \texttt{4.24e-14} & \texttt{1.39e-11} &  \texttt{1.29e-13} & \texttt{6.55e+00} & \texttt{6.55e+00} & \texttt{3.86e+02}& \texttt{5.33e-01} & \texttt{2.81e-03} \\ \hline
 $1\times 10^5$& \texttt{4.02e-13} & \texttt{9.37e-14} &  \texttt{2.04e-12} & \texttt{2.81e+00} & \texttt{2.81e+00} & \texttt{8.56e+02}& \texttt{1.16e+01} & \texttt{1.80e-02} \\ \hline
 $1\times 10^6$& \texttt{2.53e-12} & \texttt{3.54e-15} &  \texttt{3.55e-11} & \texttt{2.65e+00} & \texttt{2.65e+00} & \texttt{2.01e+03}& \texttt{1.86e+01} & \texttt{8.18e-02} \\ \hline
 $1\times 10^7$& \texttt{1.77e-11} & \texttt{1.61e-11} &  \texttt{2.44e-10} & \texttt{4.90e+00} & \texttt{4.90e+00} & \texttt{6.29e+03}& \texttt{9.44e+00} & \texttt{9.51e-01} \\ \hline
 \end{tabular} }
\end{table}


\setlength\tabcolsep{.95mm}
\begin{table}[!h]
\caption{ Experiment 3: $\mathbf{B}$ is positive semidefinite and singular with 
$[\mathbf{g}_{\parallel}]_i = 0 $ for $ 1 \le i \le r $ and 
$\| \Lambda^{\dagger} \mathbf{g}_{\parallel} \|_2 > \delta $.} {
\centering
\begin{tabular}{|c|c|c|c|c|c|c|c|c|}
   \hline $n$   	& \texttt{opt 1} &  \texttt{opt 2} &
\texttt{opt 3} & \texttt{min($\lambda(B+C_\parallel)$)} & 
$ \sigma^*_{\parallel}$ & $\sigma^*_{\perp} $	& $ \gamma $ &  time  \\
\hline
$1\times 10^3$& \texttt{1.38e-14} & \texttt{1.35e-09} &  \texttt{1.21e-14} & \texttt{1.99e+00} & \texttt{1.99e+00} & \texttt{1.45e+02}& \texttt{2.80e+00} & \texttt{3.84e-03} \\ \hline
  $ 1\times 10^4 $& \texttt{7.38e-14} & \texttt{2.98e-17} &  \texttt{4.35e-13} & \texttt{8.60e+00} & \texttt{8.60e+00} & \texttt{3.80e+03}& \texttt{1.29e+01} & \texttt{2.03e-03} \\ \hline
 $1\times 10^5$& \texttt{1.73e-13} & \texttt{8.84e-17} &  \texttt{4.17e-12} & \texttt{3.19e+00} & \texttt{3.19e+00} & \texttt{3.19e+03}& \texttt{4.67e+00} & \texttt{6.31e-03} \\ \hline
 $1\times 10^6$& \texttt{2.04e-12} & \texttt{1.22e-11} &  \texttt{4.25e-11} & \texttt{8.57e+00} & \texttt{8.57e+00} & \texttt{2.97e+04}& \texttt{1.28e+01} & \texttt{7.37e-02} \\ \hline
 $1\times 10^7$& \texttt{3.98e-11} & \texttt{7.53e-11} &  \texttt{2.42e-10} & \texttt{4.47e+00} & \texttt{4.47e+00} & \texttt{2.25e+04}& \texttt{6.63e+00} & \texttt{9.42e-01} \\ \hline\end{tabular} }
\end{table}

\setlength\tabcolsep{0.95mm}
\begin{table}[!h]
\caption{ Experiment 4: $\mathbf{B}$ is indefinite and 
$[\mathbf{g}_{\parallel}]_i = 0 $ for $ 1 \le i \le r $ with
$ \| ( \Lambda-\lambda_1\mathbf{I} )^{\dagger}\mathbf{g}_{\parallel} \|_2 > \delta $.}  {
\centering
\begin{tabular}{|c|c|c|c|c|c|c|c|c|}
   \hline $n$   	& \texttt{opt 1} &  \texttt{opt 2} &
\texttt{opt 3} & \texttt{min($\lambda(B+C_\parallel)$)} & 
$ \sigma^*_{\parallel}$ & $\sigma^*_{\perp} $	& $ \gamma $ &  time  \\
\hline
$1\times 10^3$& \texttt{1.95e-14} & \texttt{2.57e-16} &  \texttt{0.00e+00} & \texttt{2.34e+00} & \texttt{3.09e+00} & \texttt{2.38e+02}& \texttt{3.04e+00} & \texttt{3.03e-03} \\ \hline
  $ 1\times 10^4 $& \texttt{8.69e-14} & \texttt{2.16e-16} &  \texttt{0.00e+00} & \texttt{2.18e+00} & \texttt{2.59e+00} & \texttt{4.63e+02}& \texttt{2.91e+00} & \texttt{6.16e-03} \\ \hline
 $1\times 10^5$& \texttt{2.52e-13} & \texttt{4.65e-17} &  \texttt{1.72e-12} & \texttt{1.33e+01} & \texttt{1.34e+01} & \texttt{2.15e+04}& \texttt{1.98e+01} & \texttt{6.44e-03} \\ \hline
 $1\times 10^6$& \texttt{4.45e-12} & \texttt{1.24e-12} &  \texttt{1.91e-11} & \texttt{7.02e+00} & \texttt{7.21e+00} & \texttt{2.58e+04}& \texttt{1.04e+01} & \texttt{6.93e-02} \\ \hline
 $1\times 10^7$& \texttt{2.52e-11} & \texttt{5.27e-10} &  \texttt{7.46e-11} & \texttt{1.02e+00} & \texttt{1.21e+00} & \texttt{1.71e+04}& \texttt{8.35e-01} & \texttt{9.23e-01} \\ \hline \end{tabular} }
\end{table}

Tables I-VI show the results of the experiments.  In all tables,
the residual of the two optimality conditions \texttt{opt 1},
\texttt{opt 2},
and \texttt{opt 3} are on the order of $1\times10^{-10}$ or smaller.
Columns 4 in all tables show that ($\mathbf{B}+\mathbf{C}_\parallel$)
are postiive semidefinite.
Columns 6 and 7 in all the tables show that 
$\sigma_{\parallel}^*$ and  $\sigma_\perp^*$
are nonnegative.  
Thus, the solutions obtained by \SCSR{}
for these experiments satisfy the
optimality conditions to high accuracy.

Also reported in each table are the number of Newton iterations.  In the
first five experiments no more than four Newton iterations were required to
obtain $\sigma_\parallel$ to high accuracy (Column 8).  In the hard case,
no Newton iterations are required since $\sigma_\parallel^*=-\lambda_1$.
This is reflected in Table VI, where Column 4 shows that
$\sigma_\parallel^*=-\lambda_1$ and Column 8 reports no Newton iterations.)

The final column reports the time required by \SCSR{} to solve each subproblem.
Consistent with the best limited-memory methods, as $n$ gets large, the time required to solve each
  subproblem appears to grow linearly with $n$, as predicted in Section~\ref{subsec:cc}.

Additional experiments were run with
    $\mathbf{g}_\parallel\rightarrow 0$.  In particular, the experiments
    were rerun with $\bfm{g}$ scaled by factors of 
   $10^{-2}, 10^{-4},$ $10^{-6}$,
    $10^{-8}$, and $10^{-10}$.  All experiments resulted in tables similar
    to those in Tables I-VI: the optimality conditions were satisfied to
    high accuracy, no more than three Newton iterations were required in
    any experiment to find $\sigma_\parallel^*$, and the CPU times are similar
    to those found in the tables.

\setlength\tabcolsep{0.95mm}
\begin{table}[!h]
\caption{ Experiment 5: $\mathbf{B}$ is indefinite and 
$ [ \mathbf{g}_{\parallel}
]_i \ne 0$ for some $ 1 \le i \le r $.} {
\centering
\begin{tabular}{|c|c|c|c|c|c|c|c|c|}
   \hline $n$   	& \texttt{opt 1} &  \texttt{opt 2} &
\texttt{opt 3} & \texttt{min($\lambda(B+C_\parallel)$)} & 
$ \sigma^*_{\parallel}$ & $\sigma^*_{\perp} $	& $ \gamma $ &  time  \\
\hline
$1\times 10^3$& \texttt{9.11e-15} & \texttt{5.14e-16} &  \texttt{4.35e-15} & \texttt{7.54e-01} & \texttt{1.16e+00} & \texttt{1.31e+01}& \texttt{2.27e+01} & \texttt{5.60e-03} \\ \hline
  $ 1\times 10^4 $& \texttt{6.04e-14} & \texttt{8.71e-12} &  \texttt{1.25e-13} & \texttt{1.88e+00} & \texttt{2.23e+00} & \texttt{1.41e+02}& \texttt{4.15e+00} & \texttt{1.75e-03} \\ \hline
 $1\times 10^5$& \texttt{3.16e-13} & \texttt{3.27e-11} &  \texttt{2.36e-12} & \texttt{6.23e-01} & \texttt{1.24e+00} & \texttt{3.86e+02}& \texttt{4.89e+00} & \texttt{7.91e-03} \\ \hline
 $1\times 10^6$& \texttt{1.19e-12} & \texttt{0.00e+00} &  \texttt{2.82e-11} & \texttt{3.01e+00} & \texttt{3.59e+00} & \texttt{1.89e+03}& \texttt{1.77e+01} & \texttt{7.00e-02} \\ \hline
 $1\times 10^7$& \texttt{5.25e-11} & \texttt{1.02e-14} &  \texttt{1.30e-10} & \texttt{7.37e-01} & \texttt{1.43e+00} & \texttt{4.32e+03}& \texttt{1.48e+01} & \texttt{9.40e-01} \\ \hline\end{tabular} }
\end{table}

\setlength\tabcolsep{0.95mm}
\begin{table}[!h]
  \caption{ Experiment 6: $\mathbf{B}$ is indefinite and 
$[\mathbf{g}_{\parallel}]_i = 0 $ for $ 1 \le i \le r $ with
$\| \mathbf{v}_{\parallel} (
    -\lambda_1 )\|_2 \le \delta $ (the ``hard case'').} {
\centering
\begin{tabular}{|c|c|c|c|c|c|c|c|c|}
   \hline $n$   	& \texttt{opt 1} &  \texttt{opt 2} &
\texttt{opt 3} & \texttt{min($\lambda(B+C_\parallel)$)} & 
$ \sigma^*_{\parallel}$ & $\sigma^*_{\perp} $	& $ \gamma $ &  time  \\
\hline
$1\times 10^3$& \texttt{1.58e-14} & \texttt{1.21e-17} &  \texttt{2.83e-14} & \texttt{0.00e+00} & \texttt{1.09e-01} & \texttt{1.45e+02}& \texttt{1.19e+00} & \texttt{2.06e-03} \\ \hline
  $ 1\times 10^4 $& \texttt{9.07e-14} & \texttt{2.65e-17} &  \texttt{2.62e-13} & \texttt{0.00e+00} & \texttt{3.19e-01} & \texttt{4.49e+02}& \texttt{9.14e+00} & \texttt{1.31e-03} \\ \hline
 $1\times 10^5$& \texttt{8.34e-13} & \texttt{8.80e-17} &  \texttt{1.86e-12} & \texttt{0.00e+00} & \texttt{1.67e-01} & \texttt{1.45e+03}& \texttt{5.04e+00} & \texttt{4.45e-03} \\ \hline
 $1\times 10^6$& \texttt{3.87e-12} & \texttt{7.21e-17} &  \texttt{5.46e-12} & \texttt{0.00e+00} & \texttt{1.30e-01} & \texttt{3.51e+03}& \texttt{3.31e+00} & \texttt{6.77e-02} \\ \hline
 $1\times 10^7$& \texttt{4.19e-11} & \texttt{1.30e-17} &  \texttt{3.05e-10} & \texttt{0.00e+00} & \texttt{2.68e-02} & \texttt{2.81e+04}& \texttt{1.19e+01} & \texttt{9.45e-01} \\ \hline \end{tabular} }
\end{table}


\subsection{$(\mathbf{P},\infty)$-norm results}
The \SCSR{} method was tested
on randomly-generated problems of size $n=10^3$ to $n=10^7$,
organized as five experiments that test the cases enumerated 
in Algorithm~\ref{alg-pinfty}.  Since Algorithm~\ref{alg-pinfty} proceeds
componentwise (i.e., the components of $\mathbf{g}_\parallel$ and
$\mathbf{\Lambda}$
determine how the algorithm proceeds),
the experiments were designed to ensure at least one randomly-chosen
component 
satisfied the conditions of the given experiment. 
The five experiments are below:

\begin{enumerate}
\item[(E1)]  {$\left|[\mathbf{g}_\parallel]_i\right|<\delta\left|[\mathbf{\Lambda}]_{ii}\right|$ and
$[\mathbf{\Lambda}]_{ii}>\tau$.} 
 \item[(E2)] {$\left|[\mathbf{g}_\parallel]_i\right|<\tau$ and $\left|[\mathbf{\Lambda}]_{ii}\right|<\tau$.} 
\item[(E3)] {$\left|[\mathbf{g}_\parallel]_i\right|>\tau$ and $\left|[\mathbf{\Lambda}]_{ii}\right|<\tau$.} 

 \item[(E4)]{$\left|[\mathbf{g}_\parallel]_i\right|<\tau$ and $[\mathbf{\Lambda}]_{ii}<-\tau$.} 
\item[(E5)]  {$\left|[\mathbf{g}_\parallel]_i\right|>\delta\left|[\mathbf{\Lambda}]_{ii}\right|$ and
$\|\mathbf{\Lambda}]_{ii}\|>\tau$.} 
 \end{enumerate}

For these experiments, $\mathbf{S}$, $\mathbf{Y}$, and $\mathbf{g}$ were
randomly generated and then altered to satisfy the requirements described
above by each experiment.  In (E2)-(E4), $\delta$ was chosen as a
random number (in the other experiments, it was set in accordance with
the experiments' rules).  All randomly-generated vectors and matrices were
formed using the \MATLAB{} \texttt{randn} command, which draws from the
standard normal distribution.  The initial \SR{} matrix was set to
$\mathbf{B}_0=\gamma \mathbf{I}$, where $\gamma=|10*\texttt{randn(1)}|$.
Finally, the number of limited-memory updates $m$ was set to 5,
 and for simplicity, the
randomly-chosen $i$ (that defines [E1]-[E5])
was chosen to be an integer in the range $[1\,\, 5]$.

$1\times 10^3$
\begin{table}[h!]
\setlength\tabcolsep{1.5mm}
 \caption{  Results using the $(\mathbf{P},\infty)$ norm.} {
  \centering
\begin{tabular}{|l|c|c|c|}
   \hline   &  $n$   	& $ \gamma $ &  time  \\ 
\hline
  Experiment 1
            &     $1\times 10^3$& \texttt{8.76e+00} & \texttt{1.34e-03} \\ \hline
 &$ 1\times 10^4$& \texttt{1.80e-01} & \texttt{1.21e-03} \\ \hline
  &$1\times 10^5$& \texttt{7.39e+00} & \texttt{6.71e-03} \\ \hline
&$1\times 10^6$& \texttt{2.13e-02} & \texttt{1.12e-01} \\ \hline
 &$1\times 10^7$& \texttt{1.11e+01} & \texttt{1.51e+00} \\ \hline
Experiment 2 
&$1\times 10^3$& \texttt{4.47e+00} & \texttt{1.05e-03} \\ \hline
 & $ 1\times 10^4 $& \texttt{6.38e+00} & \texttt{8.74e-04} \\ \hline
 &$1\times 10^5$& \texttt{1.10e+00} & \texttt{7.37e-03} \\ \hline
 &$1\times 10^6$& \texttt{2.74e+00} & \texttt{7.94e-02} \\ \hline
 &$1\times 10^7$& \texttt{8.30e-01} & \texttt{1.39e+00} \\ \hline
Experiment 3 
&$1\times 10^3$& \texttt{2.09e+01} & \texttt{1.07e-03} \\ \hline
 & $ 1\times 10^4 $& \texttt{4.67e+00} & \texttt{9.63e-04} \\ \hline
 &$1\times 10^5$& \texttt{1.39e+01} & \texttt{6.63e-03} \\ \hline
 &$1\times 10^6$& \texttt{1.76e+01} & \texttt{7.38e-02} \\ \hline
 &$1\times 10^7$& \texttt{1.51e+01} & \texttt{1.45e+00} \\ \hline
Experiment 4 
&$1\times 10^3$& \texttt{1.08e+01} & \texttt{1.43e-03} \\ \hline
 & $ 1\times 10^4 $& \texttt{1.34e+01} & \texttt{1.06e-03} \\ \hline
 &$1\times 10^5$& \texttt{7.43e+00} & \texttt{1.23e-02} \\ \hline
 &$1\times 10^6$& \texttt{3.16e+00} & \texttt{9.00e-02} \\ \hline
 &$1\times 10^7$& \texttt{2.22e+00} & \texttt{1.41e+00} \\ \hline
Experiment 5 
&$1\times 10^3$& \texttt{1.04e+01} & \texttt{1.15e-03} \\ \hline
 & $ 1\times 10^4 $& \texttt{1.74e+01} & \texttt{9.40e-04} \\ \hline
 &$1\times 10^5$& \texttt{4.38e+00} & \texttt{1.15e-02} \\ \hline
 &$1\times 10^6$& \texttt{5.21e+00} & \texttt{9.05e-02} \\ \hline
 &$1\times 10^7$& \texttt{2.01e+00} & \texttt{1.40e+00} \\ \hline
\end{tabular} }
\end{table}

Table VII displays the results of the five experiments.  
Each experiment was run five times; the results of the third iteration
are stored in Table VII.  In all cases, the results of the third
iteration were representative of all the iterations.  The first column
of the table denotes the experiment, the second column displays the
size of the problem, and the third column reports the value of
$\gamma$.  Finally, the forth column reports the time taken to obtain the solution.

\subsection{Trust-Region Algorithm}
\label{sec:TRalg}
	\jb{In this experiment, we embed the TR subproblem solvers in a trust-region algorithm to solve unconstrained optimization problems.}
In particular, we implemented our subproblem solvers in an algorithm that is based on \cite[Algorithm 6.2]{NocW06}. A feature
of this algorithm is that the L-SR1 matrix is updated by every pair $\{(\mathbf{s}_i,\mathbf{y}_i)\}^k_{\jbo{i=k-m+1}} $
as long as $ |\mathbf{s}^T_i(\mathbf{y}_i - \mathbf{B}_{0}\mathbf{s}_i)| \ge \jbo{ \| \mathbf{s}_i \|_2 \| \mathbf{y}_i - \mathbf{B}_{i}\mathbf{s}_i \|_2}  \epsilon_{\text{SR1}} $ 
(updates are skipped if this condition is not met). \jbo{In case a full memory strategy is used (i.e, $m=\infty$)} then a \jbo{SR-1} matrix  
is updated by almost every pair $\{(\mathbf{s}_i,\mathbf{y}_i)\} $ in order to help achieve the superlinear 
convergence rate of quasi-Newton methods, \jb{in contrast to updating the matrix only when a step is accepted}.  An outline of our trust-region method implementation (Algorithm 5)
is included in the Appendix. 
\jb{In our comparisons we use the following 4 algorithms to solve the TR subproblems:
\begin{center}
\begin{tabular}{l l}
TR:SC-INF 	& Algorithm \ref{alg-pinfty} \\
TR:SC-L2 	& Algorithm \ref{alg-p2} \\
TR:L2 		& $\ell_2$-norm \cite[Algorithm 1]{BruEM15} \\ 
tr:CG 		& truncated CG \cite[Algorithm 7.5.1]{ConGT00a} \\
\end{tabular}
\end{center}
Initially, we included a $5^{\textnormal{th}}$ algorithm, LSTRS \cite{RSS08},
which performed markedly inferior to any of the above solvers and is thus not reported
as part of the outcomes in this section. We also found that Init. 2 performed significantly better than
Init. 1, and therefore report the outcomes with Init. 2 below.
Because the limited-memory updating mechanism is different whether a constant or non-constant
initialization strategy is used, we describe our results separately for C Init. and Init. 2. As part of our
comparisons we first select the best algorithm using only C Init. and only using Init. 2. Subsequently,
we compare the best algorithms to each other. In order to find default parameters for our 
best algorithms, Figures \ref{fig:default_C}, \ref{fig:default_I2_mIN}, and \ref{fig:default_I2_qIN} report results for a considerable
range of $ m$ and $q$ values.}

\jb{All remaining experiments are for the general unconstrained minimization problem}
\begin{equation}
\label{eq:min}
	\underset{ \mathbf{x} \in \mathbb{R}^n }{ \text{ minimize }  } f(\mathbf{x}),
\end{equation}
where $ f: \mathbb{R}^n \to \mathbb{R} $. We consider this problem solved
once $ \| \nabla f (\mathbf{x}_k)  \|_{\infty} \le \varepsilon $. Our convergence tolerance is set to 
be $ \varepsilon = 5\times 10^{-4} $. 
 \jbo{With $\gamma$ fixed, a L-SR1 algorithm can be implemented by only storing the matrices $ \mathbf{\Psi}_k $ and $\mathbf{M}_k^{-1}$.
In particular, with a fixed $ \gamma = \gamma_k $ in \eqref{eqn-PsiM} then $ \mathbf{M}_k^{-1} \mathbf{e}_k = \mathbf{\Psi}^T_k \mathbf{s}_k $,
so that updating the symmetric matrix $ \mathbf{M}_k^{-1} $ only uses $O(nm)$ multiplications. In this way, the overall
computational complexity and memory requirements of the L-SR1 method are reduced as compared to non-constant initializations.} 
\jb{However, using a non-constant initialization strategy can adaptively incorporate additional information, which can
be advantageous. Therefore, we compare the best algorithms for constant and non-constant
initialization strategies in Sections \ref{subsec:compC}, \ref{subsec:compI} and \ref{subsec:compAll}. }
Parameters in Algorithm 5 are
set as follows: $c_1 = 9\times 10^{-4}$, $c_2 = 0.75$, $c_3 = 0.8 $,
$c_4=2$, $c_5=0.1$, $c_6=0.75$, $c_7=0.5$ and $ \varepsilon_{\text{SR1}} = 1 \times 10^{-8} $. 

\jb{Extended performance profiles as in \cite{MahajanLeyfferKirches11} are provided. These profiles are an extension
of the well known profiles of Dolan and Mor\'{e} \cite{DolanMore02}. 
 We compare total computational time for each solver on the test set of problems.
 The performance metric $ \rho_s(\tau) $ with a given number of test problems $ n_p $ is
\begin{equation*}
		\rho_s(\tau) = \frac{\text{card}\left\{ p : \pi_{p,s} \le \tau \right\}}{n_p} \quad \text{and} \quad \pi_{p,s} = \frac{t_{p,s}}{ \underset{1\le i \le S, i \ne s}{\text{ min } t_{p,i}} },
\end{equation*} 
where $ t_{p,s}$ is the ``output'' (i.e., time) of
``solver'' $s$ on problem $p$. Here $ S $ denotes the total number of solvers for a given comparison. This metric measures
the proportion of how close a given solver is to the best result. The extended performance profiles are the same as the classical ones for $ \tau \ge 1 $.
In the profiles we include a dashed vertical grey line, to indicate $ \tau = 1 $. The solvers are compared on 62 large-scale CUTEst problems,
which are the same problems as in \cite{Burdakov2016}. Additionally, Appendix \ref{subsec:EX_QuadRos} includes supplementary comparisons on quadratics and
the Rosenbrock objectives.}

\subsection{\jb{Comparisons with constant initialization strategy (C Init.)}}
\label{subsec:compC}
\begin{figure*}[t!]	
		\includegraphics[width=\textwidth]{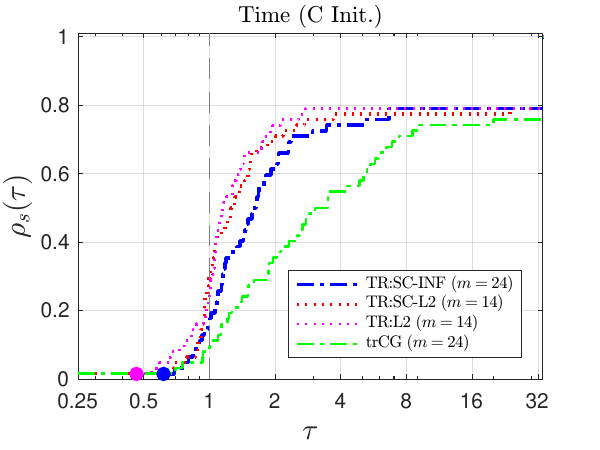}
		\caption{\jb{Comparison of best algorithms with C Init. (constant initialization), which are selected
		from Figure \ref{fig:default_C}. Observe that TR:L2 obtains the best results in this comparison. The limited-memory
		parameter $m$ is relatively large for all solvers, however since a constant initialization is used larger memory values 
		are permissible. }} 
		\label{fig:compC}       
\end{figure*}

\jb{This experiment compares the algorithms when the constant initialization C Init. from \eqref{eq:init} is used. Because
the memory allocation is essentially halved (relative to a non-constant initialization) the memory parameter $m$ includes larger values, too
(such as $m=24$). For each individual solver we first determine its optimal $ m $ parameter in Figure \ref{fig:default_C}. After selecting the
best parameters, these best solvers are then compared in Figure \ref{fig:compC}.}

\subsection{\jb{Comparisons with non-constant initialization strategy (Init. 2)}}
\label{subsec:compI}
\begin{figure*}[t!]	
		\includegraphics[width=\textwidth]{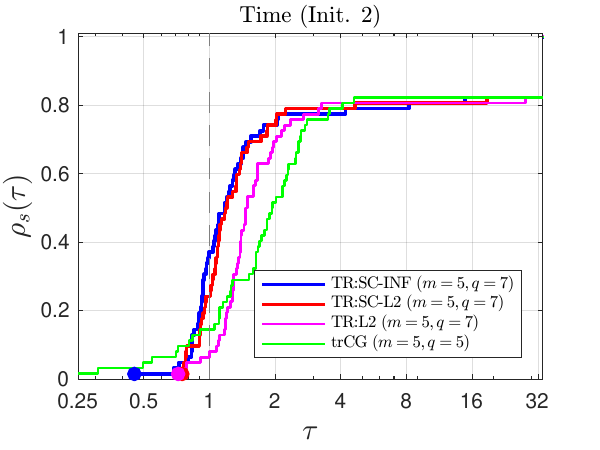}
		\caption{\jb{Comparison of best algorithms with Init.2 (non-constant initialization), which are selected
		as the best ones from Figures \ref{fig:default_I2_mIN}, and \ref{fig:default_I2_qIN}. Observe that
		TR:SC-INF and TR:SC-L2 obtain the overall best results. All algorithms use a small memory parameter $m=5$.
		Since Init. 2 is a non-constant initialization these algorithms store $ \bk{S}, \bk{Y} $ to implicitly represent 
		$ \bsk{\Psi} $, and thus the memory allocations scale with $ 2\cdot m $.}}				
		\label{fig:compI}       
\end{figure*}

\jb{Since Init. 2 depends on the parameter $q$, Figures \ref{fig:default_I2_mIN} and \ref{fig:default_I2_qIN} test each
 algorithm on a combination of $m$ and $q$ values. A comparison of the best values for each algorithm is in Figure \ref{fig:compI}.}

\subsection{\jb{Comparisons of best outcomes}}
\label{subsec:compAll}
\begin{figure*}[t!]	
		\includegraphics[width=\textwidth]{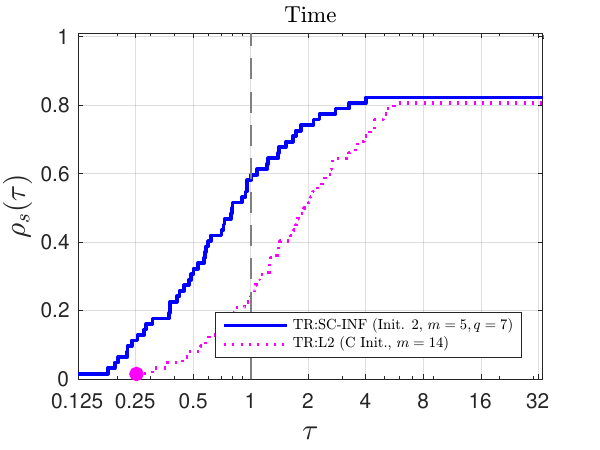}
		\caption{\jb{Overall comparison of best algorithms by selecting winners in Figures \ref{fig:compC} (C Init.) 
		and \ref{fig:compI} (Init. 2). Observe that TR:SC-INF with non-constant initialization strategy outperforms the best
		algorithm with a constant initialization (TR:L2). In sum, the trust-region algorithm with the proposed shape-changing infinity subproblem 
		solver (TR:SC-INF) obtains the best results among the comparisons on 62 large-scale CUTEst problems. }}		
		\label{fig:compALL}       
\end{figure*}

\jb{The overall best algorithms from Figures \ref{fig:compC} and \ref{fig:compI} are compared in Figure \ref{fig:compALL}.
This declares that the best performing algorithm over the sequence of experiments is TR:SC-INF with the indicated parameter values.}

\section{Concluding remarks}
In this paper, we presented a high-accuracy trust-region subproblem solver for when the Hessian is approximated by \LSR{}
matrices.  The method makes use of special shape-changing norms that decouple the original subproblem into two
separate subproblems.  Numerical experiments using the $(\mathbf{P},2)$ norm verify that solutions are computed to high accuracy in
cases when there are no closed-form solutions and also in the
so-called ``hard case''. \jb{Experiments on large-scale unconstrained optimization problems demonstrate that the proposed algorithms
perform well when compared to widely used methods, such as truncated CG or an $\ell_2$ TR subproblem algorithm.}

\appendix
\section*{APPENDIX}
\setcounter{section}{1}
\label{app:LSR1alg}
This appendix lists our implementation of the L-SR1 trust-region algorithm
from the numerical experiments in Section 5.3. This trust-region algorithm uses 
the trust-region radius adjustments from \cite[Algorithm 6.2]{NocW06} and the subproblem solvers in Algorithms 3 and 4,
as well as the orthonormal basis method (OBS) from \cite{BruEM15}. \\

\hrule 
\vspace{0.1cm}
\noindent \small{\textbf{ALGORITHM 5:} L-SR1 Shape-Changing Trust-Region Algorithms (LSR1\_SC)} 
\vspace{0.1cm}
\hrule
\begin{algorithmic}[1]
\ENSURE [$\mathbf{x}_k,\mathbf{g}_k,f_k,\text{out}$]=\texttt{LSR1\_SC}($\mathbf{x}$, $f(\mathbf{x})$, $\nabla f(\mathbf{x})$, \text{pars})
\STATE{Set constants from \text{pars}: $ 0 < c_1 < 1\times 10^{-3} $, $0 < c_2$, $ 0 < c_3 < 1 $,
$ 1 < c_4 $, $ 0 < c_5 \le c_2 $, $ 0 < c_6 < 1 $, $0 < c_7 < 1$,$ 0 < \varepsilon $, $ 0 < m $, $  0 < q $, $0 <  \varepsilon_{\text{SR1}} $,
$\text{ALG} \gets \text{pars.whichSub}$,
$\text{INIT} \gets \text{pars.whichInit}$,
$\text{SAVE} \gets \text{pars.storePsiPsi}$, ;}  
\STATE{Initialize $k \gets 0$, $k_m \gets 0$, $\mathbf{x}_k \gets \mathbf{x}$, $0<\gamma_k$, $0<\gamma_{\text{max}}$, $\text{inv}\mathbf{M}_k\gets[]$,
$\text{mIdx}\gets1:m$, $\text{iEx}\gets 0$;} 
\STATE{$f_k \gets f(\mathbf{x}_k)$, $\mathbf{g}_k \gets \nabla f(\mathbf{x}_k)$;}
\STATE{$[\mathbf{x}_{k+1},\mathbf{g}_{k+1},f_{k+1}] \gets \text{lineSearch}(\mathbf{x}_{k},\mathbf{g}_{k},f_{k})$;}
\STATE{$ \mathbf{s}_k \gets \mathbf{x}_{k+1}-\mathbf{x}_k, \mathbf{y}_k \gets \mathbf{g}_{k+1} - \mathbf{g}_k $;}
\IF {$\text{INIT}=\text{C.Init.}$}
	\STATE{\texttt{\% Constant initialization}}
	\STATE{$\gamma_k \gets \text{max}(\text{min}(\| \bfm{y}_0 \|^2 / \bfm{s}^T_0 \bfm{y}_0, \gamma_{\text{max}}),1) $}
	\STATE{$\mathbf{\Psi}_k \gets[]$;}
\ELSE
	\STATE{\texttt{\% Non-constant initialization.}} 
	\STATE{$ \gamma_k \gets \| \bfm{y}_k \|^2 / \bfm{s}^T_k \bfm{y}_k $;}
	\STATE{$\mathbf{S}_k \gets[]$,
	$\mathbf{Y}_k \gets[]$,
	$\mathbf{D}_k \gets[]$,
	$\mathbf{L}_k \gets[]$,
	$\mathbf{T}_k \gets[]$,
	$\mathbf{SS}_k \gets[]$,
	$\mathbf{YY}_k \gets[]$; }
\ENDIF
\IF {$\text{SAVE}=1$}
	\STATE{$ \bs{\Psi}\bs{\Psi}_k \gets[] $;}
\ENDIF
\STATE{$ b_k \gets  \mathbf{s}^T_k(\mathbf{y}_k - \gamma_k \mathbf{s}_k) $;}
\IF { $\varepsilon_{\text{SR1}} \jbo{\| \mathbf{s}_k \|_2  \| \mathbf{y}_k - \gamma_k \mathbf{s}_k \|_2} < \text{abs}(b_k)$}
\STATE{$k_m \gets k_m+1$;}
\STATE{$\text{inv}\mathbf{M}_k(k_m,k_m) \gets b_k$;}
\IF {$\text{INIT}=\text{C.Init.}$}
	\STATE{$[\bs{\Psi}_k,\text{mIdx}] = \texttt{colUpdate}(\bsk{\Psi},\mathbf{y}_k - \gamma_k \mathbf{s}_k,\text{mIdx},m,k) $ \texttt{\% From Procedure 1}}
\ELSE
	\STATE{$[\bk{Y},\sim] = \texttt{colUpdate}(\bk{Y},\mathbf{y}_k ,\text{mIdx},m,k) $;}
	\STATE{$[\bk{S},\text{mIdx}] = \texttt{colUpdate}(\bk{S},\mathbf{s}_k,\text{mIdx},m,k) $;}
	\STATE{$\bk{D}(k_m,k_m) = \bk{s}^T\bk{y}$;}
	\STATE{$\bk{L}(k_m,k_m) = \bk{s}^T\bk{y}$;}
	\STATE{$\bk{T}(k_m,k_m) = \bk{s}^T\bk{y}$;}
	\STATE{$\bk{SS}(k_m,k_m) = \bk{s}^T\bk{s}$;}
	\STATE{$\bk{YY}(k_m,k_m) = \bk{y}^T\bk{y}$;}
\ENDIF

\ENDIF
\STATE{$\delta_k \gets 2 \| \mathbf{s}_k \|$;}
\STATE{$k\gets k+1$;}
\WHILE {($\varepsilon \le \| \mathbf{g}_k \|_2$) \text{and} ($ k \le \text{maxIt} $) } 
\STATE{Choose TR subproblem solver to compute $ \bk{s} $ (E.g., Alg. 3, Alg. 4, $\ell_2$-norm, truncated CG);}
\STATE{\texttt{\% For example: sc\_sr1\_infty with $ \Psi^T\Psi $ updating}}
\STATE{ $\mathbf{s}_k \gets $ \texttt{sc\_sr1\_infty}($\mathbf{g}_k$,$\mathbf{S}_k(:,\text{mIdx}(1:k_m))$, 
			$\mathbf{Y}_k(:,\text{mIdx}(1:k_m))$, $\gamma_k$, $\delta_k$, $1$, $0$, $\ldots$\\
			$ \bsk{\Psi\Psi}(1:k_m,1:k_m) $,$ \text{inv}\mathbf{M}_k(1:k_m,1:k_m) $); }

\STATE{$\widehat{\mathbf{x}}_{k+1} \gets \mathbf{x}_k + \mathbf{s}_k, \widehat{f}_{k+1} \gets f(\widehat{\mathbf{x}}_{k+1}),
 \widehat{\mathbf{g}}_{k+1} \gets \nabla f(\widehat{\mathbf{x}}_{k+1})$;}
 \IF {$\text{INIT}=\text{C.Init}$}
	 \STATE{$\mathbf{b}_k(1:k_m) \gets \mathbf{\Psi}_k(:,\text{mIdx}(1:k_m))^T \mathbf{s}_k$;}
 \ELSE
 	\STATE{\texttt{\% Non-constant initialization, stores additionally $ \bk{b1}, \bk{b2} $}}
	\STATE{$\mathbf{b1}_k(1:k_m) \gets \mathbf{Y}_k(:,\text{mIdx}(1:k_m))^T \mathbf{s}_k$;}
	\STATE{$\mathbf{b2}_k(1:k_m) \gets \mathbf{S}_k(:,\text{mIdx}(1:k_m))^T \mathbf{s}_k$;}
	\STATE{$\mathbf{b}_k(1:k_m) \gets \bk{b1}(1:k_m) - \gamma_k \bk{b2}(1:k_m) $;}
 \ENDIF
 
  \STATE{$(sBs)_k \gets  \gamma_k\mathbf{s}^T_k \mathbf{s}_k + \frac{1}{2} \mathbf{b}_k(1:k_m)^T(\text{inv}
  \mathbf{M}_k(1:k_m,1:k_m) \backslash \mathbf{b}_k(1:k_m) ) $;}
  \IF{$\textnormal{INIT} = \textnormal{Init. 2}$}
  	\STATE{\texttt{\% Other non-constant initialization strategies can be implemented here}}
	\STATE{$ \gamma_k \gets \text{max}(\| \bfm{y}_{k-q} \|^2/\bfm{s}^T_{k-q} \bfm{y}_{k-q}, \cdots, \| \bk{y} \|^2/ \bk{y}^T \bk{s} )$}
  \ENDIF
  
 \STATE{$\rho_k \gets \frac{\widehat{f}_{k+1} - f_k}{\mathbf{s}^T_k \mathbf{g}_k +(sBs)_k} $; }
 \IF {$c_1 < \rho_k$}
	\STATE{$\mathbf{x}_{k+1} \gets \widehat{\mathbf{x}}_{k+1}$;}
	\STATE{$\mathbf{g}_{k+1} \gets \widehat{\mathbf{g}}_{k+1}$;}
	\STATE{$f_{k+1} \gets \widehat{f}_{k+1}$;}
\ELSE
	\STATE{$\mathbf{x}_{k+1} \gets \mathbf{x}_{k}$;}
 \ENDIF
 \IF {$c_2 < \rho_k$}
 	\IF{$\| \mathbf{s}_k \|_2 \le c_3 \delta_k$}
		\STATE{$\delta_k \gets \delta_k$;}
	\ELSE
		\STATE{$\delta_k \gets c_4 \delta_k$;}
	\ENDIF
\ELSIF {$c_5 \le \rho_k \le c_6$}
	\STATE{$\delta_k \gets \delta_k$;}
\ELSE
	\STATE{$\delta_k \gets c_7\delta_k$;}
\ENDIF
\STATE{$\mathbf{y}_k \gets \widehat{\mathbf{g}}_{k+1} - \mathbf{g}_k$;}
\STATE{$b_k \gets \mathbf{s}_k^T\mathbf{y}_k + (sBs)_k$;}
\IF{$ \varepsilon_{\text{SR1}} \jbo{\| \mathbf{s}_k \|_2  \| \mathbf{y}_k - \gamma_k \mathbf{s}_k \|_2} \le \text{abs}(b_k) $} 
	
	\IF{$\text{INIT}=\text{C.Init.}$}
		\STATE{$[\bs{\Psi}_k,\text{mIdx}] = \texttt{colUpdate}(\bsk{\Psi},\mathbf{y}_k - \gamma_k \mathbf{s}_k,\text{mIdx},m,k) $;}
	\IF{($k_m < m$)}
		\STATE{$k_m \gets k_m+1$;}
	\ENDIF
			\STATE{$\text{inv}\mathbf{M}_k(1:(k_m-1),k_m) \gets \mathbf{b}_k(1:(k_m-1))$;}
			\STATE{$\text{inv}\mathbf{M}_k(k_m,1:(k_m-1)) \gets \mathbf{b}_k(1:(k_m-1))$;}
			\STATE{$\text{inv}\mathbf{M}_k(k_m,k_m) \gets b_k$;}
		\IF {$ \textnormal{SAVE} = 1$}
			\STATE{\texttt{\% Update and store the product $ \Psi_k^T \Psi_k $}}
			\STATE{$\bsk{\Psi\Psi}(1:k_m,1:k_m) = \texttt{prodUpdate}(\bsk{\Psi\Psi},\bsk{\Psi},\bsk{\Psi},\bk{y}-\gamma_k\bk{s},\bk{y}-\gamma_k\bk{s},\text{mIdx},m,k) $;}
		\ENDIF
	\ELSE \STATE{\texttt{\% Non-constant initialization}} 
		\STATE{$[\bk{Y},\sim] = \texttt{colUpdate}(\bk{Y},\mathbf{y}_k ,\text{mIdx},m,k) $;}
		\STATE{$[\bk{S},\text{mIdx}] = \texttt{colUpdate}(\bk{S},\mathbf{s}_k,\text{mIdx},m,k) $;}
		\STATE{$\bk{T} = \texttt{prodUpdate}(\bk{T},\bk{S},0,\bk{s},\bk{y},\text{mIdx},m,k) $;}
		\STATE{$\bk{YY} = \texttt{prodUpdate}(\bk{YY},\bk{Y},\bk{Y},\bk{y},\bk{y},\text{mIdx},m,k) $;}
		
		
		
			\IF{($k_m < m$)}
				\STATE{$k_m \gets k_m+1$;}
			\ENDIF
				\STATE{$\bk{D}(k_m,k_m) \gets \bk{s}^T \bk{y}$;}		
				\STATE{$\bk{L}(k_m,1:(k_m-1)) \gets \mathbf{b1}_k(1:(k_m-1))$;}
				\STATE{$\bk{SS}(1:(k_m-1),k_m) \gets \mathbf{b2}_k(1:(k_m-1))$;}
				\STATE{$\bk{SS}(k_m,1:(k_m-1)) \gets \mathbf{b2}_k(1:(k_m-1))$;}
				\STATE{$\bk{SS}(k_m,k_m) \gets \bk{s}^T\bk{s}$;}
				\STATE{$\text{inv}\mathbf{M}_k(1:k_m,1:k_m) \gets \bk{D}(1:k_m,1:k_m)
					+ \bk{L}(1:k_m,1:k_m) + \bk{L}(1:k_m,1:k_m)^T - \gamma_k \bk{SS}(1:k_m,1:k_m)$;}
			\IF {$ \textnormal{SAVE} = 1$}
			\STATE{\texttt{\% Update and store the product $ \Psi_k^T \Psi_k $ with non-constant initialization}}
			\STATE{$\bsk{\Psi\Psi}(1:k_m,1:k_m) = \bk{YY}(1:k_m,1:k_m) - 
			\gamma_k ( \bk{T}(1:k_m,1:k_m) + \bk{T}(1:k_m,1:k_m)^T + \bk{L}(1:k_m,1:k_m) + \bk{L}(1:k_m,1:k_m)^T )
			+ \gamma_k^2 \bk{SS}(1:k_m,1:k_m) $;}
			\ENDIF
	\ENDIF 
\ENDIF 
\STATE{$\mathbf{x}_k \gets \mathbf{x}_{k+1}$, $ \mathbf{g}_{k} \gets \mathbf{g}_{k+1} $, $f_{k} \gets f_{k+1}$, $k \gets k+1$; }
\ENDWHILE
\STATE{$\text{out.numiter}\gets k, \text{out.ng} \gets \| \mathbf{g}_k \|$;}
\RETURN $\mathbf{x}_k, \mathbf{g}_k, f_k, \text{out}$
\end{algorithmic}
\hrule
 
 \subsection{\jb{Experiments to determine default parameters with constant Initialization (C Init.)}}
 \label{subsec:EXC_DEFAULT}
\begin{figure*}[t!]
	\begin{minipage}{0.48\textwidth}
		\includegraphics[trim=0 0 20 0,clip,width=\textwidth]{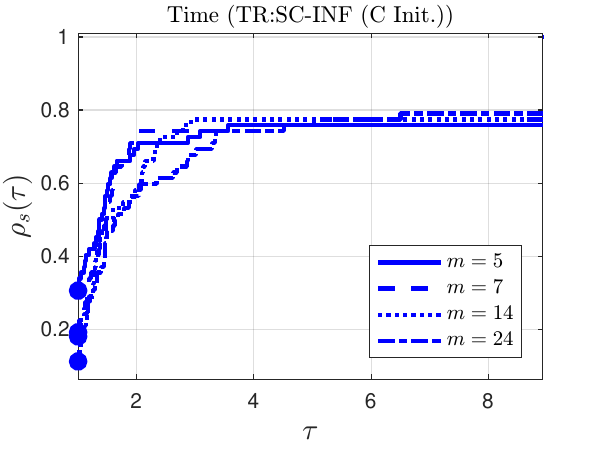}
	\end{minipage}
		\hfill
	\begin{minipage}{0.48\textwidth}
		\includegraphics[trim=0 0 20 0,clip,width=\textwidth]{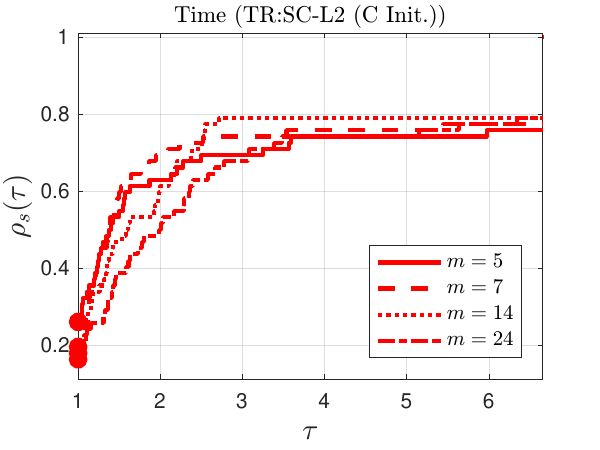}
	\end{minipage}
	\begin{minipage}{0.48\textwidth}
		\includegraphics[trim=0 0 20 0,clip,width=\textwidth]{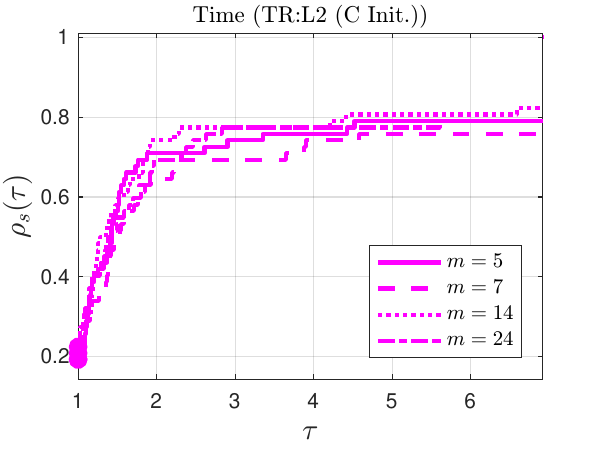}
	\end{minipage}
		\hfill
	\begin{minipage}{0.48\textwidth}
		\includegraphics[trim=0 0 20 0,clip,width=\textwidth]{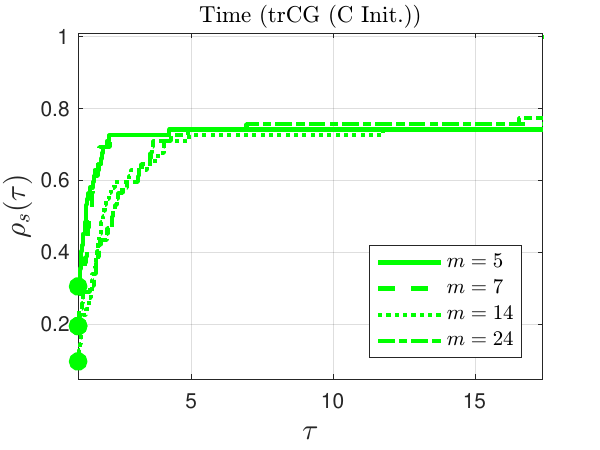}
	\end{minipage}
	\caption{\jb{Comparison of the computational times for the 4 algorithms $\{\textnormal{TR:SC-INF, TR:SC-L2, TR:L2, trCG} \}$ when a constant initialization 
	(C Init.) is used, and the limited memory parameter is $ m = [5,7,14,24] $.}} 
		\label{fig:default_C}       
\end{figure*}

 \subsection{\jb{Experiments to determine default parameters with non-constant Initialization (Init. 2)}}
 \label{subsec:EXI2_DEFAULT}

\begin{figure*}[t!]
	\begin{minipage}{0.30\textwidth}
		\includegraphics[trim=0 0 20 0,clip,width=\textwidth]{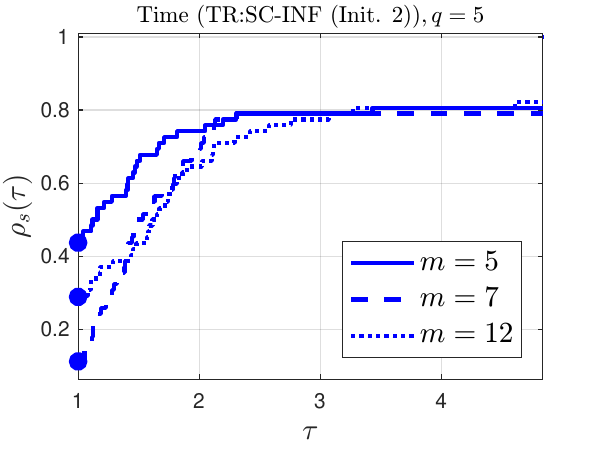}
	\end{minipage}
	\begin{minipage}{0.30\textwidth}
		\includegraphics[trim=0 0 20 0,clip,width=\textwidth]{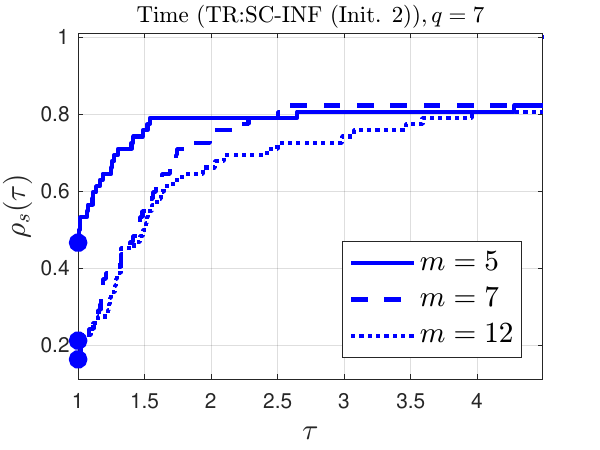}
	\end{minipage}
		\hfill
	\begin{minipage}{0.30\textwidth}
		\includegraphics[trim=0 0 20 0,clip,width=\textwidth]{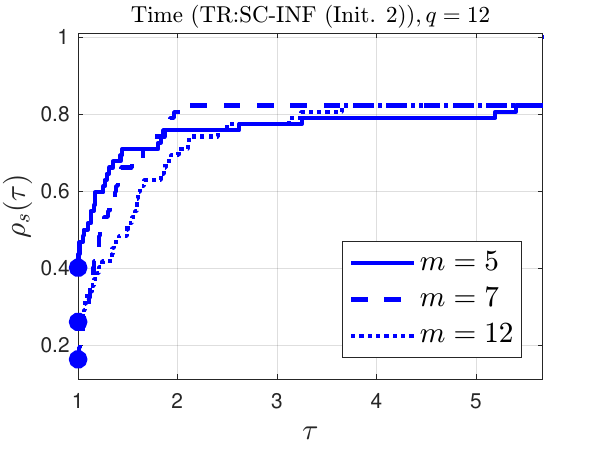}
	\end{minipage}
	\begin{minipage}{0.30\textwidth}
		\includegraphics[trim=0 0 20 0,clip,width=\textwidth]{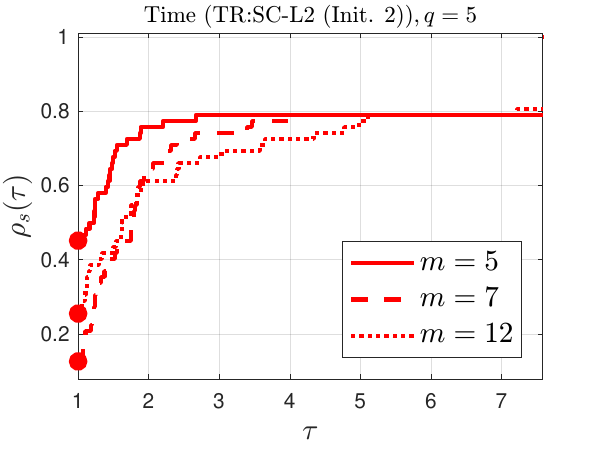}
	\end{minipage}
	\begin{minipage}{0.30\textwidth}
		\includegraphics[trim=0 0 20 0,clip,width=\textwidth]{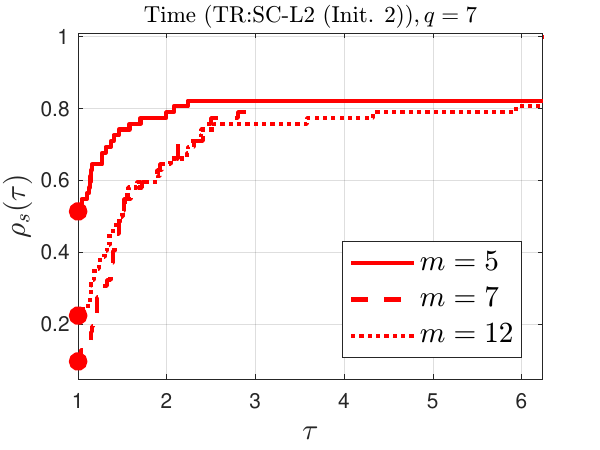}
	\end{minipage}
		\hfill
	\begin{minipage}{0.30\textwidth}
		\includegraphics[trim=0 0 20 0,clip,width=\textwidth]{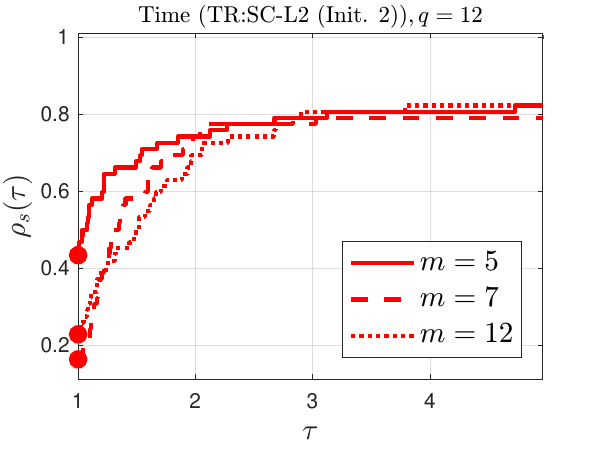}
	\end{minipage}
	\begin{minipage}{0.30\textwidth}
		\includegraphics[trim=0 0 20 0,clip,width=\textwidth]{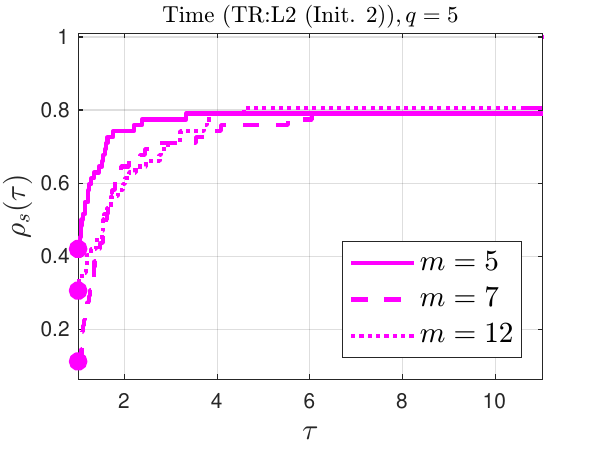}
	\end{minipage}
	\begin{minipage}{0.30\textwidth}
		\includegraphics[trim=0 0 20 0,clip,width=\textwidth]{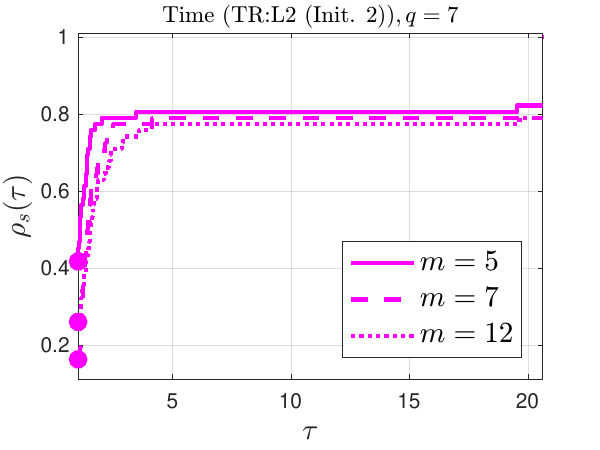}
	\end{minipage}
		\hfill
	\begin{minipage}{0.30\textwidth}
		\includegraphics[trim=0 0 20 0,clip,width=\textwidth]{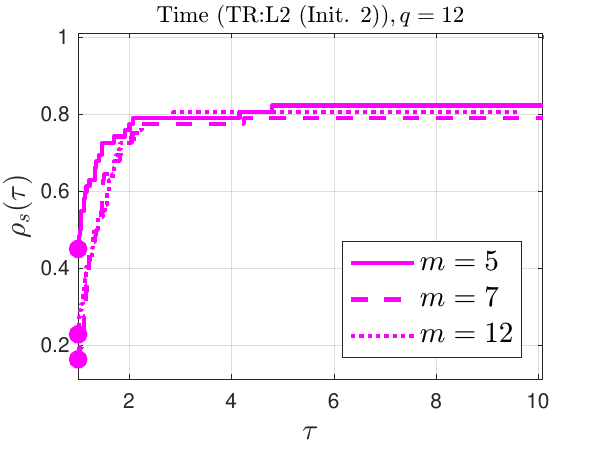}
	\end{minipage}
	\begin{minipage}{0.30\textwidth}
		\includegraphics[trim=0 0 20 0,clip,width=\textwidth]{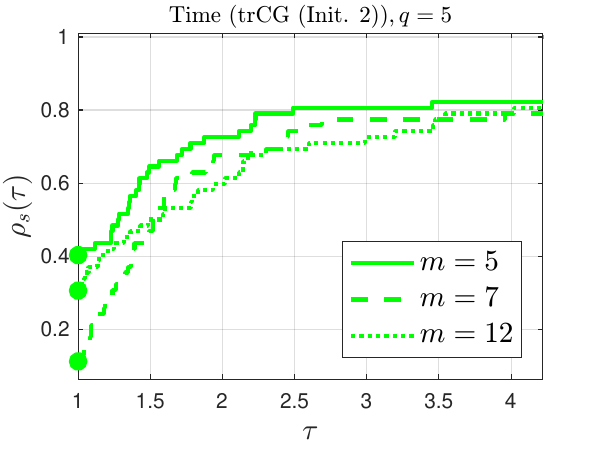}
	\end{minipage}
	\begin{minipage}{0.30\textwidth}
		\includegraphics[trim=0 0 20 0,clip,width=\textwidth]{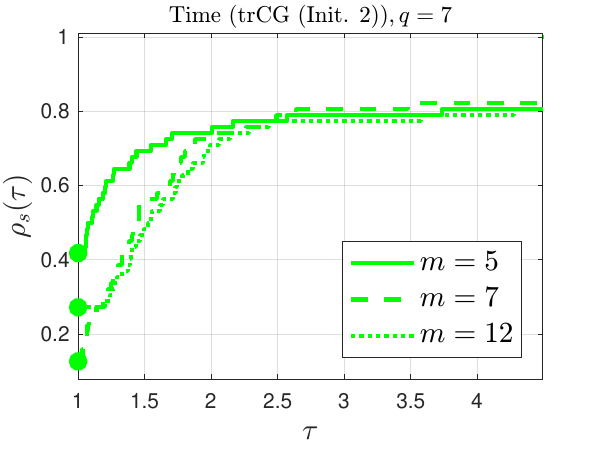}
	\end{minipage}
		\hfill
	\begin{minipage}{0.30\textwidth}
		\includegraphics[trim=0 0 20 0,clip,width=\textwidth]{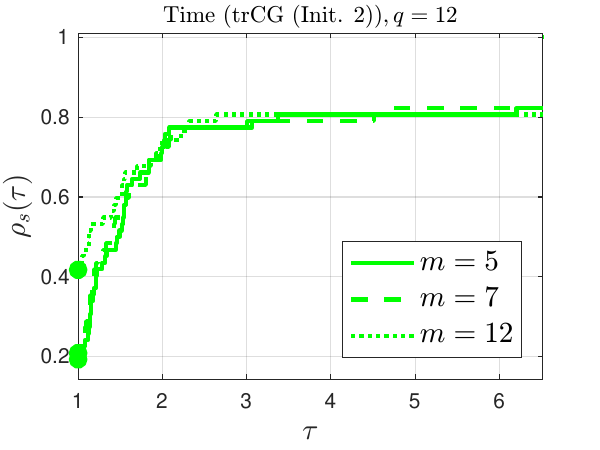}
	\end{minipage}
	\caption{\jb{Comparison of the computational times for the 4 algorithms $\{\textnormal{TR:SC-INF, TR:SC-L2, TR:L2, trCG} \}$ when the non-constant initialization 
	(Init. 2) is used, and the parameters are $q=[5,7,12]$ and $ m = [5,7,12] $.}} 
		\label{fig:default_I2_mIN}       
\end{figure*}

\begin{figure*}[t!]
	\begin{minipage}{0.30\textwidth}
		\includegraphics[trim=0 0 20 0,clip,width=\textwidth]{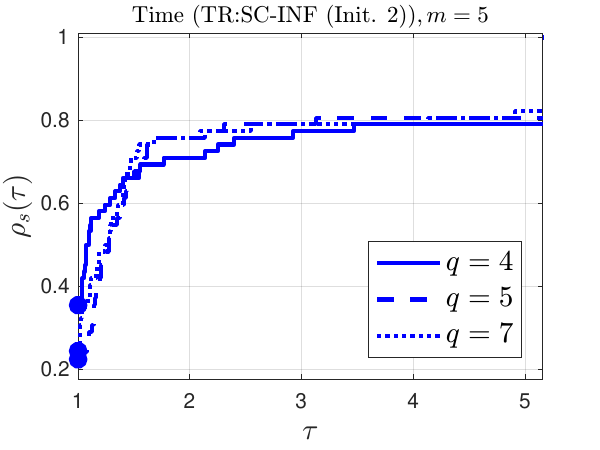}
	\end{minipage}
	\begin{minipage}{0.30\textwidth}
		\includegraphics[trim=0 0 20 0,clip,width=\textwidth]{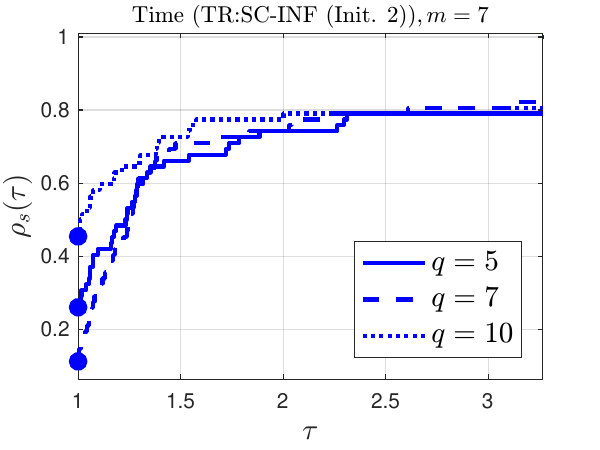}
	\end{minipage}
		\hfill
	\begin{minipage}{0.30\textwidth}
		\includegraphics[trim=0 0 20 0,clip,width=\textwidth]{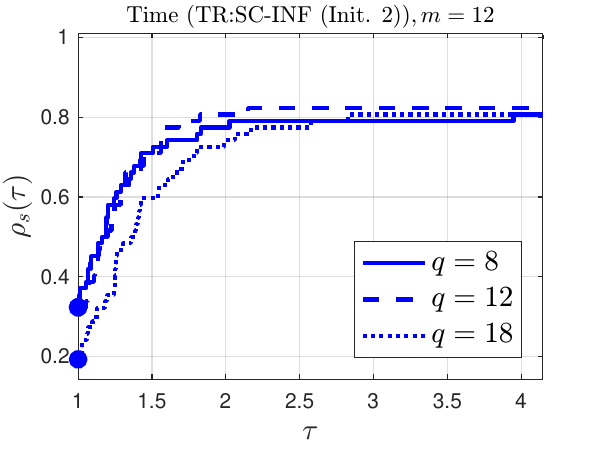}
	\end{minipage}
	\begin{minipage}{0.30\textwidth}
		\includegraphics[trim=0 0 20 0,clip,width=\textwidth]{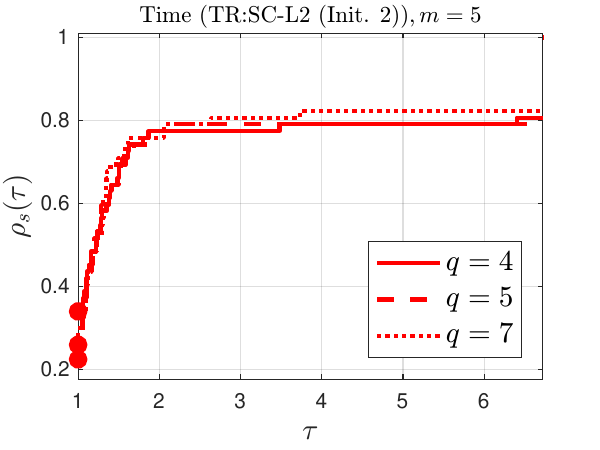}
	\end{minipage}
	\begin{minipage}{0.30\textwidth}
		\includegraphics[trim=0 0 20 0,clip,width=\textwidth]{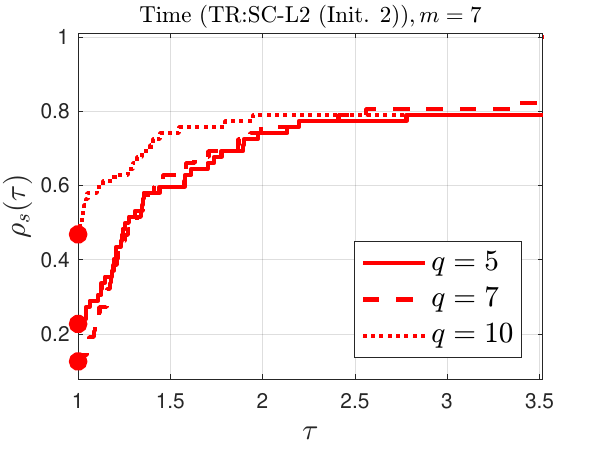}
	\end{minipage}
		\hfill
	\begin{minipage}{0.30\textwidth}
		\includegraphics[trim=0 0 20 0,clip,width=\textwidth]{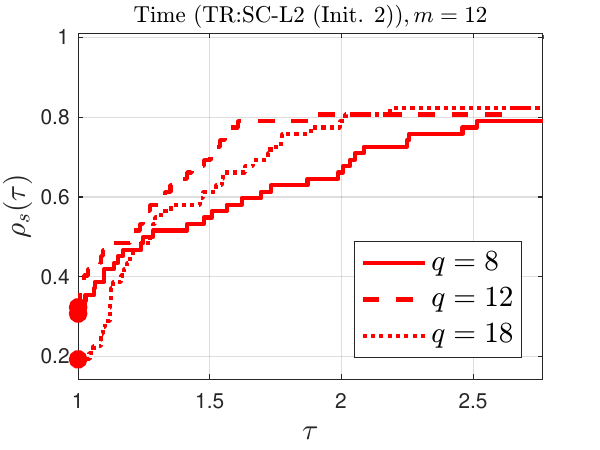}
	\end{minipage}
	\begin{minipage}{0.30\textwidth}
		\includegraphics[trim=0 0 20 0,clip,width=\textwidth]{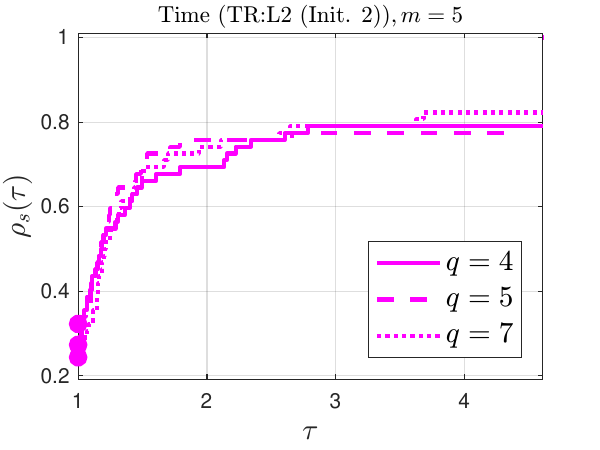}
	\end{minipage}
	\begin{minipage}{0.30\textwidth}
		\includegraphics[trim=0 0 20 0,clip,width=\textwidth]{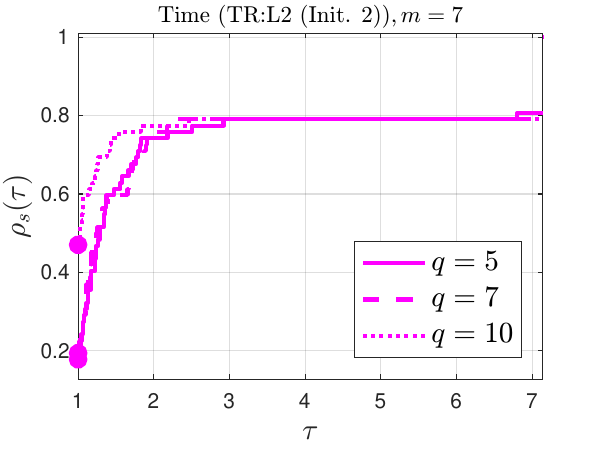}
	\end{minipage}
		\hfill
	\begin{minipage}{0.30\textwidth}
		\includegraphics[trim=0 0 20 0,clip,width=\textwidth]{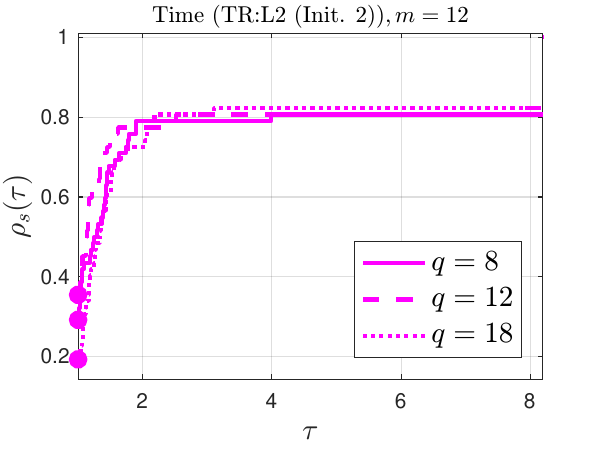}
	\end{minipage}
	\begin{minipage}{0.30\textwidth}
		\includegraphics[trim=0 0 20 0,clip,width=\textwidth]{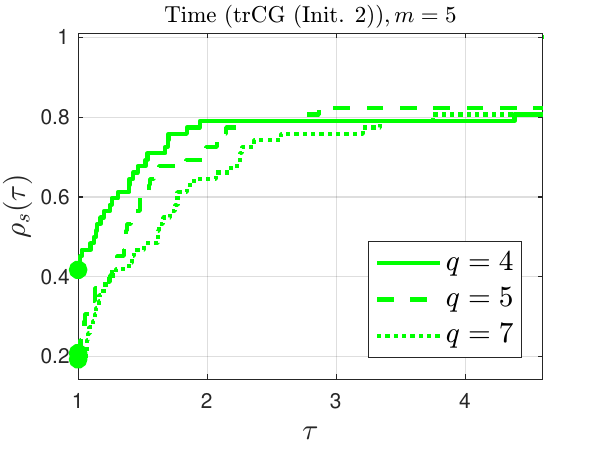}
	\end{minipage}
	\begin{minipage}{0.30\textwidth}
		\includegraphics[trim=0 0 20 0,clip,width=\textwidth]{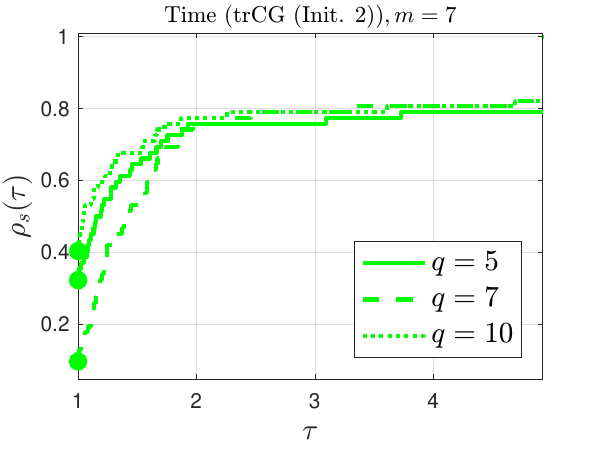}
	\end{minipage}
		\hfill
	\begin{minipage}{0.30\textwidth}
		\includegraphics[trim=0 0 20 0,clip,width=\textwidth]{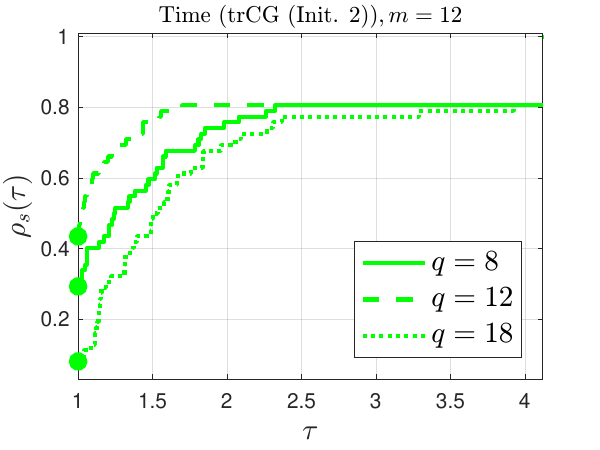}
	\end{minipage}
	\caption{\jb{Comparison of the computational times for the 4 algorithms $\{\textnormal{TR:SC-INF, TR:SC-L2, TR:L2, trCG} \}$ when the non-constant initialization 
	(Init. 2) is used, and the parameters are $m=[5,7,12]$ and $ q= [\texttt{ceil}(2/3\cdot m_i),m_i,\texttt{floor}(3/2 \cdot m_i)], 1 \le i \le 3 $.}} 
		\label{fig:default_I2_qIN}       
\end{figure*}

\subsection{\jb{Experiments on quadratics and the Rosenbrock functions}}
\label{subsec:EX_QuadRos}
In this set of experiments we vary the problem dimension as 
$ n = [5\times10^2,1\times10^3,5\times10^3,1\times10^4,5\times10^4,1\times10^5,3\times10^5] $,
set the memory parameter $m=5$, use Init. 2 for all solvers and set the maximum iterations as $ \text{maxIt} = 500 $.
In Table VIII, we let $ f (\mathbf{x}) $
be the Rosenbrock function defined by $ f(\mathbf{x}) = \sum_{i=1}^n (\mathbf{x}_{2i} - \mathbf{x}^2_{2i-1})^2 + (1-\mathbf{x}^2_{2i-1})^2  $.
 We initialize the trust-region algorithm (Algorithm 5) from the starting point $ [\mathbf{x}_0]_1 = 30, [\mathbf{x}_0]_{2:n} = 0 $. (With this initial point the gradient norm $ \| \nabla f(\mathbf{x}_0) \|_2 \approx 10^5 $). 
Table VIII reports the outcomes of using the trust-region algorithm with these three different subproblem solvers. 

\begin{table}[h!]
\label{tbl:trResults}
\setlength\tabcolsep{0.75mm} 
 \caption{  Results of solving problem \eqref{eq:min} with the Rosenbrock objective function. The maximum
 number of iterations are: $ \text{maxIt} = 500 $ and the convergence tolerance is
 $ \| \nabla f(\mathbf{x}_k) \|_{\infty} \le 1\times 10^{-4} $. The memory parameter is $ m = 5 $. The column 
 $\text{nA}$ denotes the number of ``accepted" search directions, which corresponds to line 55 in Algorithm 5 being true. Observe that
all algorithms converged to the prescribed tolerances on all problem instances.} {
  \centering
\begin{tabular}{|c| c c c c| c c c c| c c c c|}
		\hline $n$ 
		& \multicolumn{4}{c|}{TR:SC-INF (Alg. 3)} 
		& \multicolumn{4}{c|}{TR:SC-L2 (Alg. 4)}
		& \multicolumn{4}{c|}{TR-L2} \\ 
		\cline{2-13}
		& $k$ & nA & Time & $ \| \nabla f(\mathbf{x}_k) \| $ 
		& $k$ & nA & Time & $ \| \nabla f(\mathbf{x}_k) \| $ 
		& $k$ & nA & Time & $ \| \nabla f(\mathbf{x}_k) \| $ \\
\hline $5 \times 10^2$ &		\texttt{40} & \texttt{26} &  \texttt{2.00e-02} & \texttt{7.18e-05} & \texttt{46} & \texttt{29} &  \texttt{1.59e-02} & \texttt{3.11e-05} & \texttt{36} & \texttt{24} &  \texttt{1.34e-02} & \texttt{2.89e-06} \\ 
 $1 \times 10^3$ &						\texttt{38} & \texttt{24} &  \texttt{1.13e-02} & \texttt{2.23e-05} & \texttt{41} & \texttt{24} &  \texttt{1.28e-02} & \texttt{3.83e-05} & \texttt{32} & \texttt{22} &  \texttt{1.14e-02} & \texttt{2.02e-05} \\ 
 $5 \times 10^3$ &                                                \texttt{42} & \texttt{31} &  \texttt{3.02e-02} & \texttt{1.17e-05} & \texttt{38} & \texttt{29} &  \texttt{2.75e-02} & \texttt{6.38e-05} & \texttt{43} & \texttt{26} &  \texttt{4.26e-02} & \texttt{5.03e-05} \\ 
 $1 \times 10^4$ &                                                \texttt{46} & \texttt{30} &  \texttt{5.22e-02} & \texttt{4.57e-07} & \texttt{40} & \texttt{28} &  \texttt{4.20e-02} & \texttt{5.80e-05} & \texttt{48} & \texttt{29} &  \texttt{6.30e-02} & \texttt{8.87e-05} \\ 
 $5 \times 10^4$ &                                                \texttt{47} & \texttt{33} &  \texttt{2.14e-01} & \texttt{1.01e-06} & \texttt{39} & \texttt{28} &  \texttt{1.73e-01} & \texttt{1.22e-05} & \texttt{54} & \texttt{35} &  \texttt{2.85e-01} & \texttt{6.92e-05} \\ 
 $1 \times 10^5$ &                                                \texttt{40} & \texttt{31} &  \texttt{3.94e-01} & \texttt{6.82e-05} & \texttt{58} & \texttt{39} &  \texttt{4.81e-01} & \texttt{1.06e-05} & \texttt{44} & \texttt{27} &  \texttt{4.97e-01} & \texttt{1.57e-08} \\ 
 $3 \times 10^5$ &                                                \texttt{60} & \texttt{39} &  \texttt{2.74e+00} & \texttt{1.63e-06} & \texttt{53} & \texttt{33} &  \texttt{2.49e+00} & \texttt{3.52e-06} & \texttt{68} & \texttt{43} &  \texttt{3.53e+00} & \texttt{1.70e-05} \\ 
\hline
\end{tabular}} 
\end{table}

In table IX, we let $ f (\mathbf{x}) $
be quadratic functions defined by $ f(\mathbf{x}) = \mathbf{g}^T \mathbf{x} + \frac{1}{2}( \mathbf{x}^T( \phi \mathbf{I} + 
\mathbf{Q} \mathbf{D} \mathbf{Q}^T) \mathbf{x})  $. In particular, we let $ \mathbf{Q} \in \mathbb{R}^{n \times r} $
be a rectangular matrix and $ \mathbf{D} \in \mathbb{R}^{r \times r} $ be a diagonal matrix. We initialize the trust-region algorithm (Algorithm 5) from the starting point $ \mathbf{x}_0 = \mathbf{0} $. We generate $ \mathbf{Q} = \texttt{rand}(n,r) $, $ \mathbf{D} = \text{diag}(\texttt{rand}(r,1)) $ and
$ \mathbf{g} = \texttt{randn}(n,1) $, after initializing the random number generator by the 
command $ \texttt{rng}(\texttt{`default'}) $. Moreover, we set $ r =10 $, $ \phi = 100 $ and the maximum number of iterations as $ \text{maxIt} = 500 $. 
All other parameters of the method are as before. Table IX reports the outcomes of using the trust-region algorithm with the three different subproblem solvers.

\begin{table}[h!]
\label{tbl:trResults}
\setlength\tabcolsep{0.75mm} 
 \caption{  Results of solving problem \eqref{eq:min} with quadratic objective functions. The maximum
 number of iterations are set as $ \text{maxIt} = 500 $ and the convergence tolerance 
 $ \| \nabla f(\mathbf{x}_k) \|_{\infty} \le 1\times 10^{-4} $. The memory parameter is $ m = 5 $. The column 
 $\text{nA}$ denotes the number of ``accepted" search directions (line 55 in Algorithm 5 is true). Observe that
Alg. 3 and Alg. 4 converged on all problems. Moreover, Alg. 3 and Alg. 4 were fastest on the on the largest two problem instances.} {
  \centering
\begin{tabular}{|c| c c c c| c c c c| c c c c|}
		\hline $n$ 
		& \multicolumn{4}{c|}{TR:SC-INF (Alg. 3)} 
		& \multicolumn{4}{c|}{TR:SC-L2 (Alg. 4)}
		& \multicolumn{4}{c|}{TR:L2 ($\ell_2 $ \cite{BruEM15})} \\ 
		\cline{2-13}
		& $k$ & nA & Time & $ \| \nabla f(\mathbf{x}_k) \| $ 
		& $k$ & nA & Time & $ \| \nabla f(\mathbf{x}_k) \| $ 
		& $k$ & nA & Time & $ \| \nabla f(\mathbf{x}_k) \| $ \\
\hline $5 \times 10^2$ & 		\texttt{8} & \texttt{6} &  \texttt{5.75e-02} & \texttt{6.56e-06} & \texttt{8} & \texttt{6} &  \texttt{2.66e-02} & \texttt{6.56e-06} & \texttt{6} & \texttt{4} &  \texttt{2.89e-02} & \texttt{8.03e-06} \\ 
$1 \times 10^3$ &	 			 \texttt{8} & \texttt{6} &  \texttt{7.25e-03} & \texttt{4.06e-05} & \texttt{8} & \texttt{6} &  \texttt{7.51e-03} & \texttt{4.06e-05} & \texttt{6} & \texttt{4} &  \texttt{6.67e-03} & \texttt{5.11e-05} \\ 
$5 \times 10^3$ & 	 		 \texttt{21} & \texttt{15} &  \texttt{2.47e-02} & \texttt{8.96e-05} & \texttt{21} & \texttt{15} &  \texttt{2.83e-02} & \texttt{9.14e-05} & \texttt{16} & \texttt{10} &  \texttt{2.79e-02} & \texttt{3.71e-05} \\ 
$1 \times 10^4$ & 	 		 \texttt{23} & \texttt{18} &  \texttt{3.86e-02} & \texttt{7.79e-05} & \texttt{23} & \texttt{18} &  \texttt{3.65e-02} & \texttt{5.21e-05} & \texttt{19} & \texttt{14} &  \texttt{3.65e-02} & \texttt{4.28e-05} \\ 
$5 \times 10^4$ &	 			 \texttt{45} & \texttt{33} &  \texttt{2.16e-01} & \texttt{1.58e-05} & \texttt{60} & \texttt{46} &  \texttt{2.28e-01} & \texttt{9.71e-05} & \texttt{27} & \texttt{21} &  \texttt{1.20e-01} & \texttt{9.13e-05} \\ 
$1 \times 10^5$ & 	 		 \texttt{62} & \texttt{49} &  \texttt{5.04e-01} & \texttt{9.80e-05} & \texttt{79} & \texttt{64} &  \texttt{5.86e-01} & \texttt{9.72e-05} & \texttt{500} & \texttt{494} &  \texttt{4.05e+00} & \texttt{4.09e-04} \\ 
$3 \times 10^5$ & 	 		 \texttt{20} & \texttt{15} &  \texttt{8.99e-01} & \texttt{3.49e-05} & \texttt{22} & \texttt{17} &  \texttt{8.37e-01} & \texttt{3.86e-05} & \texttt{26} & \texttt{17} &  \texttt{1.11e+00} & \texttt{9.97e-05} \\ 
\hline
\end{tabular}} 
\end{table}
Remarkably, observe in the outcomes of Tables VIII and IX that a limited memory trust-region algorithm using our subproblem solvers is able to solve
large optimization problems, with $ n \approx 1 \times 10^5 $, within seconds. Moreover, we observe that the proposed algorithms (Algorithm 3 and Algorithm 4) may require 
fewer iterations on some problems than a $\ell_2$-norm method and use less computational time. Future research, can investigate
the effectiveness of a L-SR1 trust-region algorithm for non-convex objective functions and improve 
on the efficiency of the implementation.



\bibliographystyle{unsrt}
\bibliography{myrefs}
                             
\end{document}